\documentstyle[amscd, amstex]{amsart}

\newtheorem{thm}{Theorem}[section]
\newtheorem{lem}[thm]{Lemma}
\newtheorem{cor}[thm]{Corollary}
\newtheorem{pr}[thm]{Proposition}

\theoremstyle{definition}
\newtheorem{rem}[thm]{Remark}

\newtheorem{defn}[thm]{Definition}

\newcommand{\psx}{{\bf P}_{\Sigma_X}}
\newcommand{\tta}{{\bf T}_\tau}
\newcommand{\ts}{{\bf T}_\sigma}
\newcommand{\mult}{{\rm mult}}
\newcommand{\im}{{\rm im}}
\newcommand{\key}{\bibitem}
\newcommand{\psd}{{{\bf P}_{\Sigma_D}}}
\newcommand{\ps}{{{\bf P}_{\Sigma}}}
\newcommand{\vol}{{\rm vol}}
\newcommand{\gr}{{\rm Gr}}
\newcommand{\res}{{\rm Res}}
\newcommand\hidot{{\raise1pt\hbox{$\scriptscriptstyle\bullet$}}}
\newcommand\lodot{{\raise.3pt\hbox{$\scriptscriptstyle\bullet$}}}
\newcommand{\dd}{{\rm d}}
\newcommand{\pp}{{\bf P}}
\newcommand{\ttt}{{\bf T}}
\newcommand{\tr}{{\rm Tr}}
\newcommand{\inte}{{\rm int}}
\newcommand{\pic}{{\rm Pic}}
\begin{document}

\title 
[The Hodge structure of semiample hypersurfaces]
{The Hodge structure of semiample hypersurfaces and a generalization of the monomial-divisor  mirror map}
\author{Anvar R. Mavlyutov}
\address {Max-Planck-Institut f\"ur Mathematik, Bonn, D-53111, Germany.}
\email{anvar@@mpim-bonn.mpg.de}

\thanks{This research was conducted by the author for the Clay Mathematics Institute.}

\begin{abstract}
We explicitly compute 
the cohomology ring of semiample nondegenerate
hypersurfaces in complete simplicial toric
varieties. The monomial-divisor  mirror map 
is generalized to a map between the whole Picard group and the space of infinitesimal
deformations for a mirror pair of Calabi-Yau hypersurfaces.
This map is compatible with certain vanishing limiting products of the subrings
of the chiral rings,
on which  the ring structure is related to a product of the roots of $A$-type Lie
algebra.
\end{abstract}

\maketitle

\tableofcontents

In this paper we solve the long-standing problem of describing the cohomology  of semiample hypersurfaces in complete simplicial toric
varieties. The problem became important when physicists found a relation between quantum field theory and 
the geometry of Calabi-Yau manifolds. Toric varieties used as an ambient space  provide a large number of examples
of Calabi-Yau manifolds realized as complete intersections. 
Cohomology of Calabi-Yau manifolds is interpreted  as a finite dimensional subspace of the Hilbert space of states in physics.
Physicists used the (mathematically non-rigorous) orbifold Landau-Ginzburg theory  
to compute the cohomology of Calabi-Yau manifolds in some examples
coming from toric varieties (e.g., \cite{v}).
 However, no general solution in physics was found. 
About ten years ago physicists discovered Mirror Symmetry, which started with a certain symmetry of the Hodge numbers
of two topologically distinct Calabi-Yau manifolds. Soon after that, V.~Batyrev in \cite{b1}
generalized the  construction of physicists
using a duality of reflexive polytopes and Calabi-Yau hypersurfaces in toric varieties.
Being the  resolutions of ample hypersurfaces, these mirror symmetric hypersurfaces are  semiample.
The problem to compute the cohomology is also relevant to the calculation of the B-model chiral ring since
the cohomology is isomorphic to this ring. 

In most discussions of the Mirror Theorems (see \cite{ck}), 
it is not mentioned that the
results proved by Givental and Lian-Liu-Yau  
apply only to some (called toric
and polynomial) parts  of the A-model and B-model chiral rings of a mirror pair of
Calabi-Yau manifolds.  An explicit
description of the full cohomology of the manifolds is a necessary part to extend their
techniques to a complete solution.

The history of calculating the cohomology of hypersurfaces starts
with a paper of P.~Griffiths (see \cite{g}). There he found the residue (primitive)  part  of 
the middle cohomology of smooth hypersurfaces in a projective space. Such hypersurfaces are automatically
ample, implying certain vanishing theorems which allow to do the calculations.
In a later paper \cite{cg}, the product structure was computed on the residue part.
Using essentially the same method, V.~Batyrev and D.~Cox in \cite{bc} described  the residue part of
the middle cohomology of ample quasismooth hypersurfaces in complete simplicial toric varieties,
for which the analogous vanishing theorems hold as well.
For semiample hypersurfaces the vanishing theorems do not hold -- that is where the complexity of
the problem arises. However, the paper \cite{bc} also contained an alternative method 
of finding the description for 
an open subset of ample quasismooth hypersurfaces which consists of nondegenerate 
(transversal to the
torus orbits) hypersurfaces. In this case, the Gysin spectral sequence together with
the ampleness have been used.
We noticed that this spectral sequence can be applied to the semiample nondegenerate hypersurfaces as well,
but because of the missing ampleness property the complexity still appears.
While this paper provides a solution only for the semiample nondegenerate hypersurfaces,
we said in \cite{m3} that our strategy is to use this solution in order to find an answer to the question
for all semiample quasismooth hypersurfaces. The answer should be very similar since the semiample
quasismooth hypersurfaces are flat deformations of the nondegenerate ones.

Here is the plan of our paper. In Section~\ref{s:ar}, we start with introducing the notation, definitions, and then review some basic
results from \cite{m2}.

Section~\ref{s:cohreg} is devoted to deriving our main result: an explicit description of the cohomology ring
of semiample nondegenerate hypersurfaces in complete simplicial toric varieties.
In \cite{m2}  we showed that this cohomology splits into the toric and residue parts.  The toric part 
is the image of cohomology of the ambient space, while its complement  comes 
{from} the residues of rational differential forms with poles along the hypersurface.
The toric part of the cohomology has already been calculated, and
we follow the algorithmic approach suggested in \cite{m2} to derive the residue part.
In parallel, we also compute the product structure on the cohomology.

In Section~\ref{s:pic}, we apply the results from Section~\ref{s:cohreg} 
to find the Picard group of semiample nondegenerate
hypersurfaces. We give two alternative descriptions of the  Picard group, one of which is algebraic
 and the other is geometric.
A previous approach of \cite{r} to compute the  Picard group  has gaps in the proof, complicated
restrictions and not quite explicit (the toric variety ${\Bbb P}_{(\Sigma_0,L)}$ in the proof of
\cite[Theorem~2]{r} is not necessarily simplicial, hence \cite[Proposition~4]{r} could not be used
in this case).
 
In Section~\ref{s:bcor}, we use one of the main results from \cite{m2} to 
describe a subring of  $H^*(X,\wedge^*{\cal T}_X)$ of  semiample minimal Calabi-Yau 
hypersurfaces in toric varieties, which contains
$H^1(X,{\cal T}_X)$. Here, $\wedge^*{\cal T}_X$ is the {\it Zariski $p$-th  exterior power 
of the tangent sheaf}, introduced in \cite{m2}. 
We discover that 
the remarkable product structure on the subring is closely 
related to a  product of the roots of an $A$-type Lie algebra.
Looking at the subring of the cohomology of semiample minimal Calabi-Yau hypersurfaces, generated by 
the Picard group, we find that a very similar  product structure takes place.

Using the above, in Section~\ref{s:gmd},
 we propose a  generalization of the monomial-divisor  mirror map of Aspinwall, Greene and Morrison,
which is supposed to be the derivative of the mirror map at the large radius limit points.
Our map is supported by a compatibility of certain vanishing limiting products of the chiral rings
of Calabi-Yau hypersurfaces in the Batyrev mirror construction.
We also predict that the same vanishing occurs in the quantum cohomology of the Calabi-Yau hypersurfaces. 

{\it Acknowledgments.}  While all the essential techniques for solving this problem were developed in my 
Ph.D. dissertation (see \cite{m3}),
 I found the last clue this summer being a Liftoff researcher of the Clay Mathematics Institute.
I would like to thank  David Cox for suggesting me this complicated and interesting
problem for my dissertation research, and, also, for many helpful discussions and 
useful comments.

\section{A brief review of semiample divisors in toric varieties.}\label{s:ar}

In this section we  review some basic notation and facts from
\cite{m2} about  semiample divisors in toric varieties.

Throughout this paper we will use the following dictionary:\\
$M$ is a lattice of rank $d$.\\
$N=\text{Hom}(M,{\Bbb Z})$ is the dual lattice. \\
$M_{\Bbb R}$ and $N_{\Bbb R}$ are the $\Bbb R$-scalar extensions of $M$ and $N$, respectively.\\
$\Sigma$ is a finite rational (usually simplicial)   fan in $N_{\Bbb R}$.\\
${\bf P}_{\Sigma}$ is a  $d$-dimensional toric variety associated with   $\Sigma$. \\
${\bf T}_\sigma$ is a torus corresponding to the cone $\sigma\in\Sigma$.\\
$V(\sigma)$ is a toric variety equal to the closure of  ${\bf T}_\sigma$ in $\ps$.\\
$\Sigma(k)$ is the set of all $k$-dimensional cones in $\Sigma$. \\
$\Sigma(1)=\{\rho_1,\dots,\rho_n\}$ is the set of $1$-dimensional
cones in $\Sigma$ with the minimal integral generators 
$e_1,\dots,e_n$, respectively. \\
$D_i$ is a torus invariant divisor in ${\bf P}_\Sigma$, corresponding to $\rho_i$ . \\
$S=S(\ps)={\Bbb C}[x_1,\dots,x_n]$ is the homogeneous coordinate ring of $\ps$.\\
$\Delta_D=\{m\in M_{\Bbb R}:\langle m,e_i\rangle\geq-a_i
\text{ for all } i\}\subset M_{\Bbb R}$ is a convex  polytope of a Weil divisor $D=\sum_{i=1}^{n}a_iD_i$.\\

A Cartier divisor $D$ on $\ps$ is called {\it semiample} if
${\cal O}_{\ps}(D)$ is generated by global sections.
The Iitaka dimension of a Cartier divisor $D$  on $\ps$ is defined by
$$\kappa(D):=\kappa({\cal O}_{\ps}(D))=\dim\phi_D(\ps),$$
where $\phi_D:\ps@>>>{\Bbb P}(H^0(\ps,{\cal O}_{\ps}(D)))$ is the rational 
map defined by the sections of the line bundle ${\cal O}_{\ps}(D)$. 
For a torus invariant $D$, $\kappa(D)$ coincides with  the dimension of the associated polytope $\Delta_D$.

\begin{defn} \cite{m2} A semiample divisor $D$ on a complete toric variety $\ps$ 
is called {\it $i$-semiample} if  the  Iitaka dimension $\kappa(D)=i$.
\end{defn}

The following result which describes the unique properties of 
the map associated with a semiample divisor is a simplification of Theorem~1.4 in \cite{m2}.

\begin{thm}\label{t:fun}  Let $[D]\in A_{d-1}(\ps)$ be an $i$-semiample divisor class on  
a complete toric variety $\bold P_\Sigma$ of dimension $d$.
Then, there exists a unique complete toric variety $\psd$ with 
a surjective morphism
$\pi:{\bf P}_\Sigma@>>>{\bf P}_{\Sigma_D}$, arising {from} 
a surjective homomorphism of lattices $\tilde\pi:N@>>>N_D$ which maps the fan 
$\Sigma$ into $\Sigma_D$, 
such that $\pi^*[Y]=[D]$ for some ample divisor $Y$ on $\psd$.
Moreover,  $\dim\psd=i$, and the fan 
$\Sigma_D$ is the normal fan of $\Delta_D$ for a torus invariant $D$.
\end{thm}

\begin{rem} 
The fan $\Sigma_D$ lies in the space $(N_D)_{\Bbb R}:=N_{\Bbb R}/N'_{\Bbb R}$, 
where $N'=\{v\in N:\psi_D(-v)=-\psi_D(v)\}$ is a sublattice of $N$
and $\psi_D$ is the support function associated with $D$. Also,  $\Delta_D$ lies in $(M_D)_{\Bbb R}$, where $M_D:={N'}^\perp\cap M$
 (for details see \cite[Section~1]{m2}). 
\end{rem}

\begin{rem} 
The  map $\pi$ in the above theorem  is well known 
and the variety  $\psd$ can be represented as 
${\bold Proj}\bigl(\oplus_{k\ge0}H^0(\ps,{\cal O}_\ps(kD))\bigr)$.
The condition on the map  says that it is  surjective with connected fibers.
It would be interesting to see an analogue of the above theorem in a broader context.
\end{rem}

We will need to use the following intersection properties of semiample divisors.

\begin{lem}\label{l:int}
 If $D$ is an $i$-semiample divisor on a complete toric variety  
$\ps$, then the intersection number
$(D^k\cdot V(\tau))>0$ for any $\tau\in\Sigma(d-k)$, such that $\tilde\pi(\tau)$ is
contained in a cone of $\Sigma_D(i-k)$, 
and $(D^k\cdot V(\tau))=0$ for all other $\tau\in\Sigma(i-k)$.
\end{lem}

\begin{pf} The arguments are the same as for Lemma~1.4 in \cite{m1}.
\end{pf}

Let $D$ be a semiample (torus invariant) divisor in degree 
$\beta\in A_{d-1}(\ps)$ for  a complete toric variety $\ps$.
By Theorem~\ref{t:fun}, we have
the associated toric morphism $\pi:{\bf P}_\Sigma\rightarrow{\bf P}_{\Sigma_D}$.
Let also $\sigma\in\Sigma_X$ be the smallest cone , 
containing the image of the cone $\tau\in\Sigma$, and 
$\sigma'\in\Sigma(d-i+\dim\sigma)$ be such that
$\tau\subset\sigma'$ and $\tilde\pi(\sigma')\subset\sigma$.
{From} \cite[Section~1]{m2} it follows that there is a  natural commutative diagram:
$$\minCDarrowwidth0.7cm
\begin{CD}
S(\Sigma)_{p\beta} @>\varphi^*_{\sigma'}>> {S(V(\sigma'))}_{p\beta^{\sigma'}} @>{\pi_{\sigma'}}_*>> 
S(V(\sigma))_{p\beta^\sigma}\\
@VVV        @VVV             @VVV \\
 H^0(\pp_{\Sigma},{\cal O}_{\pp_{\Sigma}}(pD))\hspace{-0.1cm}@>\varphi^*_{\sigma'}>>\hspace{-0.1cm} H^0(V(\sigma'),{\cal O}_{V(\sigma')}(pD_{\sigma'}))
\hspace{-0.1cm}@>{\pi_{\sigma'}}_*>>\hspace{-0.1cm} 
H^0(V(\sigma),{\cal O}_{V(\sigma)}(p\pi_*D_{\sigma'})),
\end{CD}$$
where $\beta^{\sigma'}=[D_{\sigma'}]$, $D_{\sigma'}:=D|_{V(\sigma')}=D\cdot V(\sigma')$,
 in the Chow group of $V(\sigma')$, $\beta^\sigma=\pi_*[D_{\sigma'}]$,
$\varphi^*_{\sigma'}$ is the restriction of the global sections, the map ${\pi_{\sigma'}}_*$ is
a push-forward of the sections,
and the vertical arrows are isomorphisms.

In this paper, we are interested in the {\it nondegenerate} hypersurfaces $X\subset\ps$, which have only transversal intersections
with the tori $\ts$. This is equivalent to saying that
$x_1\partial f/\partial x_1,\dots,x_n\partial f/\partial x_n$ do not vanish simultaneously on $\ps$, where
$f\in S_\beta$ determines the hypersurface $X$. A generic hypersurface in a given semiample degree is nondegenerate, by \cite[Proposition~6.8]{d}.
We have the following property for such hypersurfaces.

\begin{pr}\label{p:regh} \cite{m2} 
Let $X$ be a semiample nondegenerate   hypersurface in a complete toric variety $\ps$, and let
$\pi:\ps@>>>\psx$ be the associated morphism  for $[X]\in A_{d-1}(\ps)$,
then  $Y=\pi(X)$ is a nondegenerate ample hypersurface, and $X=\pi^{-1}(Y)$. 
\end{pr}

For the affine hypersurfaces in tori, there is a Lefschetz-type theorem.

\begin{thm}\label{t:lef}
Let $X$ be an $i$-semiample nondegenerate  hypersurface in a complete  toric variety $\ps$.
Then the natural restrictions
$$H^l(\ttt)@>>>H^l(X\cap\ttt)$$
are isomorphisms for $l<i-1$, where $\ttt\subset\ps$ is the maximal dimensional torus.
\end{thm}

\begin{pf} This follows from  the corresponding result (see \cite{dk}), when the hypersurface is ample,
equation (2) in \cite{m1} and the K\"unneth isomorphism.
\end{pf}

In \cite[Section~5]{m2}, we proved that the cohomology of a semiample nondegenerate hypersurface $X$ in a complete
simplicial toric variety $\ps$ has a direct sum decomposition:
$$H^*(X)=H^*_{\rm toric}(X)\oplus H^*_{\rm res}(X),$$
where 
$H^*_{\rm toric}(X):=\im(H^*(\ps)@>i^*>>H^*(X))$ and 
$H^*_{\rm res}(X):=\im(H^{*+1}(\ps\setminus X)@>\res>>H^*(X))$
are called the {\it toric} and {\it residue} parts of the cohomology.
These parts are  orthogonal to each other:
$$\int_X H^*_{\rm toric}(X)\cup H^*_{\rm res}(X)=0.$$
In \cite[Theorem~5.1]{m2}, we also calculated 
$$H^*_{\rm toric}(X)\cong H^*(\ps)/Ann([X])$$
where $Ann([X])$ is the annihilator of the class $[X]\in H^2(\ps)$.
The cohomology of $\ps$  is isomorphic to 
$${\Bbb C}[D_1,\dots,D_n]/(P(\Sigma)+SR(\Sigma)),$$
where the generators correspond to the torus invariant divisors of $\ps$, 
and where
$$P(\Sigma)=\biggl\langle \sum_{i=1}^n \langle m,e_i\rangle D_i: m\in M\biggr\rangle,$$
$$SR(\Sigma)=\bigl\langle D_{i_1}\cdots D_{i_k}:\{e_{i_1},\dots,e_{i_k}\}\not\subset\sigma 
\text{ for all }\sigma\in\Sigma\bigr\rangle.$$
Hence, $H^{*}_{\rm toric}(X)$ is isomorphic to  the bigraded ring
$$A_1(X)_{*,*}:={\Bbb C}[D_1,\dots,D_n]/I,$$
where $I=(P(\Sigma)+SR(\Sigma)):[X]$ is the ideal quotient,
and $D_k$ have the degree $(1,1)$.

We finish this section with some additional frequently used notation:\\
$S(V(\tau))$, for $\tau\in\Sigma$, is the homogeneous coordinate ring of $V(\tau)$.\\
$\varphi_\tau:X\cap V(\tau)\subset X$ is the inclusion for $\tau\in\Sigma$.\\
$\varphi_{\tau,\gamma}:X\cap V(\gamma)\subset X\cap V(\tau)$ is the inclusion  for $\tau,\gamma\in\Sigma$ such that 
$\tau\subset\gamma$.\\
$f_\tau=\varphi_\tau^*(f)\in S(V(\tau))$ denotes the polynomial defining 
$X\cap V(\tau)$ for $\tau\in\Sigma$.\\
$\inte\tau=\tau\setminus\bigl(\bigcup_{\gamma\prec\tau}\gamma\bigr)$
is the relative interior of a cone $\tau$, where the union is by all proper faces of $\tau$. \\
$\pi_\tau:V(\tau)@>>>V(\sigma)$ is the restriction of $\pi$ 
for $\tau\in\Sigma$ and $\sigma\in\Sigma_X$
 such that $\tilde\pi(\inte\tau)\subset\inte\sigma$.\\
$f_\sigma={\pi_\tau}_*(f_\tau)\in S(V(\sigma))$ is the polynomial defining the hypersurface
$\pi(X)\cap V(\sigma)$ in $V(\sigma)$.\\
$d(\tau):=d-\dim\tau$ for $\tau\in\Sigma$, and $i(\sigma):=i-\dim\sigma$ for $\sigma\in\Sigma_X$.\\
$\Gamma_\sigma\subset\Delta_D$ is the face corresponding to $\sigma\in\Sigma_D$.\\
$\beta_1^\sigma=\deg(\prod_{\tilde\pi(\rho_k)\subset\sigma}x_k)\in A_{d-1}(\ps)$ for $\sigma\in\Sigma_X$.

\section{The Hodge structure of semiample nondegenerate hypersurfaces.}\label{s:cohreg}

The goal of this  section is to compute the residue part of  the cohomology of  semiample nondegenerate hypersurfaces.
This will completely generalize our previous results in \cite{m1} and \cite{m2}, concerning the middle cohomology of big and
nef hypersurfaces. 
The essential idea comes from the algorithmic approach in \cite[Section~5]{m2} to  computing of
 the residue  part of cohomology. In fact, our proof will use an induction, where \cite[Theorem~4.4]{m1} may be considered
as the base for this induction.
 This allows us to get an explicit description of the generators 
of cohomology.  
Quite remarkably, the Poincar\'e duality does the rest of the job: it  gives all the relations among 
the generators. 

Let $X$ be an $i$-semiample  nondegenerate hypersurface  in a complete simplicial toric variety $\ps$, and let 
$\pi:\ps@>>>\psx$ be the  associated contraction.
{From} the diagram (14) in \cite{m2}, we have an  exact sequence:
$$\bigoplus_{k=1}^n H_{\rm res}^{s-2}(X\cap D_k) @>\oplus{\varphi_{\rho_k}}_!>>H_{\rm res}^s(X)@>>>\gr_s^W PH^s(X\cap\ttt)@>>>0.$$
A similar sequence can be written for $X\cap D_k\subset D_k$ instead of $X\subset\ps$ as well as for all 
 hypersurfaces $X\cap V(\tau)$ in the toric subvarieties $V(\tau)$:
\begin{equation}\label{e:exa}
\bigoplus_{\tau\subset \tau'\in\Sigma(\dim\tau+1)} H_{\rm res}^{s-2}(X\cap V(\tau')) @>>>H_{\rm res}^s(X\cap V(\tau))
@>>>\gr_s^W PH^s(X\cap\tta)@>>>0.
\end{equation}

Our first step is to show how to split the residue part $H_{\rm res}^s(X\cap V(\tau))$. 
Then it will be clear how to obtain
 a description
of $H_{\rm res}^s(X)$ from the above exact sequences. We closely follow the arguments in \cite[Section~6]{m2}.
By the construction of the fan $\Sigma_X$ there is the smallest cone $\sigma\in\Sigma_X$, 
containing the image of the cone $\tau\in\Sigma$. The map $\pi$ restricted to the toric subvariety $V(\tau)\subset\ps$ gives 
the contraction $\pi_\tau:V(\tau)@>>>V(\sigma)$. 
 We will  need  to verify some  information
 with a help of a simplicial toric subvariety $V(\sigma')\subset V(\tau)$ associated with
 a cone $\sigma'\in\Sigma(d-i+\dim\sigma)$ such that
$\tau\subset\sigma'$ and the image of $\sigma'$ is in $\sigma$.
In this case, the toric variety $V(\sigma')$ is mapped birationally onto $V(\sigma)$, and
 we get a commutative diagram
\begin{equation}\label{e:pulb}
\begin{CD}
H^0(V(\sigma),\Omega^{i(\sigma)}_{V(\sigma)}(*Y_\sigma)) @>\pi_\tau^*>>
H^0(V(\tau),\Omega^{i(\sigma)}_{V(\tau)}(*X_\tau)) \\
@.\hspace{-2.5cm}\searrow^{\!\!\scriptstyle\pi_{\sigma'}^*}\hspace{2cm} @VV{\scriptstyle{\scriptstyle\varphi_{\tau,\sigma'}^*}}V \\
@. H^0(V(\sigma'),\Omega^{i(\sigma)}_{V(\sigma')}(*X_{\sigma'})),                   
\end{CD}
\end{equation}
where we use short notations: $i(\sigma)=i-\dim\sigma$, $Y_\sigma=Y\cap V(\sigma)$, $X_\tau=X\cap V(\tau)$,
$X_{\sigma'}=X\cap V(\sigma')$.
There is a natural map 
$$S(V(\sigma))_{(p+1)\beta^{\sigma}-\beta_0^{\sigma}}@>>>H^0(V(\sigma),\Omega^{i(\sigma)}_{V(\sigma)}((p+1)Y_{\sigma}))$$  
sending a polynomial $A$ to the form
$A\Omega_{V(\sigma)}/f_{\sigma}^{p+1}$, where $\Omega_{V(\sigma)}$ is the $i(\sigma)$-form for
 the complete toric variety 
$V(\sigma)$, corresponding to a basis $m_1^\sigma,\dots,m^\sigma_{i(\sigma)}$ of the lattice $M_X\cap\sigma^\perp$ (see \cite[Definition~9.3]{bc}).
The complete toric variety $V(\sigma)$ has homogeneous coordinates
$y_{\gamma}$ corresponding to the integral generators $e_{\gamma}$ of the 1-dimensional cones of the fan of $V(\sigma)$.
We want to compute the pull-back of the form $A\Omega_{V(\sigma)}/{(f_{\sigma})}^{p+1}$ with respect to the maps in (\ref{e:pulb}).
Since the form vanishes outside the torus $\ts$, 
the calculation can be done using the affine coordinates  
$$t^\sigma_j=\prod_{\gamma}y_{\gamma}^{\langle m^\sigma_j,e_{\gamma}\rangle},\quad
j=1,\dots,i(\sigma),$$ 
on the torus. 
Under the map $\pi_{\sigma'}$, the coordinates $t^\sigma_j$ are identified with 
the monomials 
$\prod_{\gamma'}x_{\gamma'}^{\langle m^\sigma_j,e_{\gamma'}\rangle}$, where $x_{\gamma'}$, for 
$\sigma'\subset\gamma'\in\Sigma(\dim(\sigma')+1)$,
 are the homogeneous coordinates on $V(\sigma')$, corresponding to the integral generators $e_{\gamma'}$.
Note
\begin{multline*}
\Omega_{V(\sigma)}=\prod_{\gamma}y_{\gamma}
\biggl(\sum_{\gamma}\langle m^\sigma_1,e_{\gamma}\rangle\frac{\dd y_{\gamma}}{y_{\gamma}}\biggr)
\wedge\cdots\wedge 
\biggl(\sum_{\gamma}\langle m^\sigma_{i(\sigma)},e_{\gamma}\rangle\frac{\dd y_{\gamma}}{y_{\gamma}}\biggr)
\\
=\prod_{\gamma}y_{\gamma}\frac{\dd t^\sigma_1}{t^\sigma_1}\wedge\cdots\wedge
\frac{\dd t^\sigma_{i(\sigma)}}{t^\sigma_{i(\sigma)}}.
\end{multline*}
Therefore, we get
\begin{multline*}
\pi_{\sigma'}^*\fracwithdelims(){A\Omega_{V(\sigma)}}{f_{\sigma}^{p+1}}=
\pi_{\sigma'}^*\fracwithdelims(){A\prod_{\gamma}y_{\gamma}}{f_{\sigma}^{p+1}}
\frac{\dd t^\sigma_1}{t^\sigma_1}\wedge\cdots\wedge
\frac{\dd t^\sigma_{i(\sigma)}}{t^\sigma_{i(\sigma)}}
\\=
\frac{\pi_{\sigma'}^*(A\prod_{\gamma}y_{\gamma})}{f_{\sigma'}^{p+1}}
\biggl(\sum_{\gamma'}\langle m^\sigma_1,e_{\gamma'}\rangle\frac{\dd x_{\gamma'}}{x_{\gamma'}}\biggr)
\wedge\cdots\wedge 
\biggl(\sum_{\gamma'}\langle m^\sigma_{i(\sigma)},e_{\gamma'}\rangle\frac{\dd x_{\gamma'}}{x_{\gamma'}}\biggr)
\\
=\frac{\pi_{\sigma'}^*(A\prod_{\gamma}y_{\gamma})}{(\prod_{\gamma'}x_{\gamma'})f_{\sigma'}^{p+1}}\Omega_{V(\sigma')},
\end{multline*}
where  $\Omega_{V(\sigma')}$ is the $i(\sigma)$-form for the complete simplicial toric variety 
$V(\sigma')$, corresponding to the basis $m_1^\sigma,\dots,m^\sigma_{i(\sigma)}$ of the lattice 
$M\cap{\sigma'}^\perp=M_X\cap\sigma^\perp$ (see \cite[Definition~9.3]{bc}).

Similarly to the above, we can also compute the pull-back  $\pi_{\tau}^*$ of
the form $A\Omega_{V(\sigma)}/f_{\sigma}^{p+1}$. But now  the coordinates
$t^\sigma_j$  can be pulled by the map $\pi_{\tau}$ to the torus $\ttt_\tau$ and expressed
in terms of the homogeneous coordinates of $\ps$:
$$t^\sigma_j=\prod_{k=1}^nx_k^{\langle m^\sigma_j,e_k\rangle}.$$
Completing $m_1^\sigma,\dots,m^\sigma_{i(\sigma)}$ to a basis of the lattice $M\cap{\tau}^\perp$, we
get a full set of coordinates  
$t^\sigma_1,\dots,t^\sigma_{i(\sigma)},t^\sigma_{i(\sigma)+1},\dots,t^\sigma_{d-\dim\tau}$ on the torus
$\ttt_\tau$. Changing the basis, if necessary, we may also assume that these coordinates are well defined on 
the other tori which map onto $\ts$. On all such tori the map  $\pi_{\tau}$ is given by projecting to the first 
$i(\sigma)$ affine coordinates. Consequently, we get 
\begin{multline*}
\pi_{\tau}^*\fracwithdelims(){A\Omega_{V(\sigma)}}{f_{\sigma}^{p+1}}=
\pi_{\tau}^*\fracwithdelims(){A\prod_{\gamma}y_{\gamma}}{f_{\sigma}^{p+1}}\frac{\dd t^\sigma_1}{t^\sigma_1}\wedge\cdots\wedge
\frac{\dd t^\sigma_{i(\sigma)}}{t^\sigma_{i(\sigma)}}
\\
=\frac{\pi_{\tau}^*(A\prod_{\gamma}y_{\gamma})}{{f}^{p+1}}
\biggl(\sum_{k=1}^n\langle m^\sigma_1,e_k\rangle\frac{\dd x_k}{x_k}\biggr)
\wedge\cdots\wedge 
\biggl(\sum_{k=1}^n\langle m^\sigma_{i(\sigma)},e_k\rangle\frac{\dd x_k}{x_k}\biggr)
\\
=\frac{\pi_{\tau}^*(A\prod_{\gamma}y_{\gamma})}{(\prod_{\tilde\pi(\rho_k)\not\subset\sigma}x_k){f}^{p+1}}
\Omega_\sigma,
\end{multline*}
where we introduce  an $i(\sigma)$-form on the  toric variety $V(\tau)$:
$$\Omega_\sigma:=\prod_{\tilde\pi(\rho_k)\not\subset\sigma}x_k
\biggl(\sum_{k=1}^n\langle m^\sigma_1,e_k\rangle\frac{\dd x_k}{x_k}\biggr)
\wedge\cdots\wedge 
\biggl(\sum_{k=1}^n\langle m^\sigma_{i(\sigma)},e_k\rangle\frac{\dd x_k}{x_k}\biggr).$$
Here, we slightly abused the notation using 
the identification
$$\bigl(S/\langle x_k:\,\tilde\pi(\rho_k)\subset\sigma\rangle\bigr)_{(p+1)\beta}
\cong S(V(\tau))_{(p+1)\beta},$$
which follows the same way as the equation (17) in \cite{m2}.

{From} the above calculations we can see a commutative diagram:
$$
\begin{CD}
S_{(p+1)\beta-\beta_0+\beta_1^\sigma}@>>>H^0(V(\tau),\Omega^{i(\sigma)}_{V(\tau)}((p+1)X_\tau))\\
@VV{\scriptstyle{\scriptstyle\varphi_{\sigma'}^*}}V  @VV{\scriptstyle{\scriptstyle\varphi_{\tau,\sigma'}^*}}V  \\
S(V(\sigma'))_{(p+1)\beta^{\sigma'}-\beta_0^{\sigma'}}@>>>H^0(V(\sigma'),\Omega^{i(\sigma)}_{V(\sigma')}((p+1)X_{\sigma'})),
\end{CD}
$$
where $\beta_1^\sigma=\deg(\prod_{\tilde\pi(\rho_k)\subset\sigma}x_k)$ and 
the horizontal arrows are represented  by the forms 
$$\frac{B\Omega_\sigma}{f^{p+1}}\,\text{ and }\,\frac{B'\Omega_{V(\sigma')}}{f_{\sigma'}^{p+1}}$$
for the corresponding polynomials $B$ and $B'$, while the arrow on the left is defined by 
$$\varphi_{\sigma'}^*(B):=\frac{\varphi_{\sigma'}^*(B\prod_{\tilde\pi(\rho_k)\not\subset\sigma}x_k)}
{\prod_{\gamma'}x_{\gamma'}}.$$
Applying the residue map, we get another diagram:
\begin{equation}\label{e:rdia}
\begin{CD}
S_{(p+1)\beta-\beta_0+\beta_1^\sigma}@>\res_\tau>>H^p(X\cap V(\tau),\Omega^{i(\sigma)-1-p}_{X\cap V(\tau)})\\
@VV{\scriptstyle{\scriptstyle\varphi_{\sigma'}^*}}V  @VV{\scriptstyle{\scriptstyle\varphi_{\tau,\sigma'}^*}}V  \\
S(V(\sigma'))_{(p+1)\beta^{\sigma'}-\beta_0^{\sigma'}}@>\res_{\sigma'}>>H^p(X\cap V(\sigma'),\Omega^{i(\sigma)-1-p}_{X\cap V(\sigma')}).
\end{CD}
\end{equation} 
By \cite[Theorem~4.4]{m1}, we know that the bottom arrow induces a well defined map on 
$R_1(f_{\sigma'})_{(p+1)\beta^{\sigma'}-\beta_0^{\sigma'}}$. 
To show a similar result for the top arrow, as in \cite[Definition~6.5]{m2}, we introduce the following rings.

\begin{defn}\label{d:rs1}
Given  $f\in S_\beta$ of semiample degree $\beta=[D]\in A_{d-1}(\ps)$   and $\sigma\in\Sigma_D$, 
let
$J^\sigma_0(f)$ be the ideal  in $S$  
generated by the ideal $J_0(f)$ and all $x_k$ such that $\tilde\pi(\rho_k)\subset\sigma$, and
let $J^\sigma_1(f)$ be the ideal quotient $J^\sigma_0(f):(\prod_{\tilde\pi(\rho_k)\not\subset\sigma}x_k)$. 
Then we get the quotient rings $R_0^\sigma(f)=S/J_0^\sigma(f)$ and  $R_1^\sigma(f)=S/J_1^\sigma(f)$
graded by the Chow group $A_{d-1}(\ps)$.
\end{defn}

\begin{pr}\label{p:iso} 
Let $\beta=[D]\in A_{d-1}(\ps)$ be $i$-semiample, and let $\sigma\in\Sigma_D$ and  $\sigma'\in\Sigma(d-i+\dim\sigma)$
be such that  $\tilde\pi(\sigma')\subset\sigma$.
If $\beta_0^{\sigma'}=\deg(\prod_{\gamma'}x_{\gamma'})\in A_{i(\sigma)-1}(V(\sigma'))$ and
$\beta_0^{\sigma}=\deg(\prod_{\gamma}y_{\gamma})\in A_{i(\sigma)-1}(V(\sigma))$ denote the anticanonical degrees, 
then, there are natural isomorphisms induced by $\varphi_{\sigma'}^*$ and ${\pi_{\sigma'}}_*$:
 
{\rm (i)} $R^\sigma_0(f)_{*\beta}\,{\cong}\, R_0(f_{\sigma'})_{*\beta^{\sigma'}}
\,{\cong}\,R_0(f_{\sigma})_{*\beta^{\sigma}}$,

{\rm (ii)} $R^\sigma_1(f)_{*\beta-\beta_0+\beta_1^\sigma}\,{\cong}\,R_1(f_{\sigma'})_{*\beta^{\sigma'}-\beta_0^{\sigma'}}
\,{\cong}\,R_1(f_{\sigma})_{*\beta^{\sigma}-\beta_0^{\sigma}}$.

{\rm (iii)} $R^\sigma_0(f)_{*\beta-\beta_0+\beta_1^\sigma}\,{\cong}\,R_0(f_{\sigma'})_{*\beta^{\sigma'}-\beta_0^{\sigma'}}
\,{\cong}\,R_0(f_{\sigma})_{*\beta^{\sigma}-\beta_0^{\sigma}}$.
\end{pr}

\begin{pf} This is just a generalization of Proposition~6.6 in \cite{m2}, and the proof follows the same way.
\end{pf}

In order to prove that the top arrow in (\ref{e:rdia}) is well defined on $R^\sigma_1(f)_{(p+1)\beta-\beta_0+\beta_1^\sigma}$
we will need  to use an explicit representation of  the $(i(\sigma)-1-p,p)$ Hodge component of the residue
$\res(B\Omega_\sigma/f^{p+1})\in H^{i(\sigma)-1}(X\cap V(\tau))$.

\begin{pr}\label{p:coc}
Let $X\subset\ps$ be an $i$-semiample nondegenerate hypersurface 
defined by $f\in S_\beta$. Given $\tau\in\Sigma$ and  the smallest cone $\sigma\in\Sigma_X$,   containing
the image of $\tau$, then,
 under the natural map
$$\check{H}^{p}({\cal U}|_{X\cap V(\tau)},\Omega_{X\cap V(\tau)}^{i(\sigma)-1-p})\rightarrow 
H^{p}(X\cap V(\tau),\Omega_{X\cap V(\tau)}^{i(\sigma)-1-p})\cong H^{i(\sigma)-1-p,p}(X\cap V(\tau)),$$   
the Hodge component of the residue
 $\res_\tau(A\Omega_\sigma/f^{p+1})=\res(A\Omega_\sigma/f^{p+1})^{i(\sigma)-1-p,p}$, for
 $A\in S_{(p+1)\beta-\beta_0+\beta_1^\sigma}$, is represented by the \v{C}ech cocycle
$$\frac{(-1)^{i(\sigma)-1+(p(p+1)/2)}}{p!}
\fracwithdelims\{\}{AK_{i_{p}}\cdots K_{i_0}\Omega_\sigma}{f_{i_0}\cdots f_{i_{p}}}_{i_0\dots i_{p}}
\in C^p({\cal U}|_{X\cap V(\tau)},\Omega^{i(\sigma)-1-p}_{X\cap V(\tau)}),$$
where ${\cal U}=\{{U}_i\}_{i=1}^n$ with ${U}_i=\{x\in\ps:x_i f_i(x)\ne0\}$, 
$f_i:=\partial f/\partial x_i$,
is the open cover of $\ps$, and where $K_i$  
is the contraction operator $\frac{\partial}{\partial x_{i}}\lrcorner$.
\end{pr}
 
\begin{pf}
This proposition generalizes our previous results 
(see \cite[Theorem~3.3]{m1} and \cite[Proposition~7.2]{m2}) as well as the corresponding results
in \cite{cg}.
However, here we suggest another way to do calculations. There is no need to plunge into hypercohomology as it was done in
\cite{cg}. Since we are interested only in the particular Hodge component of the residue, the two exact sequences 
\begin{multline*}
0@>>>\Omega_{V(\tau)}^{i(\sigma)-p}(\log X_\tau)@>>>\Omega_{V(\tau)}^{i(\sigma)-p}(X_\tau)@>\dd>>
\Omega_{V(\tau)}^{i(\sigma)-p+1}(2X_\tau)/\Omega_{V(\tau)}^{i(\sigma)-p+1}(X_\tau)@>\dd>>
\\
\cdots@>\dd>>
\Omega_{V(\tau)}^{d(\tau)}((d(\tau)-i(\sigma)+p+1)X_\tau)/
\Omega_{V(\tau)}^{d(\tau)}((d(\tau)-i(\sigma)+p)X_\tau)@>>>0
\end{multline*}
(here, $d(\tau):=\dim(V(\tau))=d-\dim\tau$; see \cite[Proposition~10.1]{bc}) and
$$0@>>>\Omega_{V(\tau)}^{i(\sigma)-p}@>>>\Omega_{V(\tau)}^{i(\sigma)-p}(\log X_\tau)@>\res>>
\Omega_{X_\tau}^{i(\sigma)-p-1}@>>>0$$
suffice to produce a representative of the component.
Define  auxiliary sheaves $L_s$ as the kernels of the morphisms
$$
\Omega_{V(\tau)}^{i(\sigma)-p+s-1}(sX_\tau)/\Omega_{V(\tau)}^{i(\sigma)-p+s-1}((s-1)X_\tau)
@>\dd>>
\Omega_{V(\tau)}^{i(\sigma)-p+s}((s+1)X_\tau)/\Omega_{V(\tau)}^{i(\sigma)-p+s}(sX_\tau)
$$
for $s=1,\dots,p+1$.
Then we get the short exact sequences
$$
0@>>>L_s@>>>\Omega_{V(\tau)}^{i(\sigma)-p+s-1}(sX_\tau)/\Omega_{V(\tau)}^{i(\sigma)-p+s-1}((s-1)X_\tau)
@>>>L_{s+1}@>>>0$$
giving rise to the maps in cohomology:
\begin{multline}\label{m:seq}
H^0(V(\tau),L_{p+1})@>>>H^1(V(\tau),L_p)@>>>
\\
\cdots@>>>H^{p-1}(V(\tau),L_2)@>>>
H^p(V(\tau),\Omega_{V(\tau)}^{i(\sigma)-p}(\log X_\tau)).
\end{multline} 

The first part of our proof is to show that the form 
$A\Omega_\sigma/f^{p+1}$ is a section in  $H^0(V(\tau),L_{p+1})$. For this it suffices to check that the differential 
of this form is a form on $V(\tau)$ with the order of the pole along $X_\tau$ at most $p+1$.
Since
$$
\dd\fracwithdelims(){A\Omega_\sigma}{f^{p+1}}=\frac{\dd(A\Omega_\sigma)}{f^{p+1}}
-\frac{(p+1)A\dd f\wedge\Omega_\sigma}{f^{p+2}},
$$
 the desired property would follow from the fact:
\begin{equation}\label{e:mul}
{\rm d}f\wedge\Omega_\sigma\equiv0\quad\mbox{ modulo multiples of } f,  x_k\text{ and }\dd x_k \text{ for }\tilde\pi(\rho_k)\subset\tau.
\end{equation} 
In order to prove this equation,  note that $\Omega_\sigma$ is equal up to a multiple by a constant to 
$$\frac{x_{l_1}\cdots x_{l_{\dim(\sigma')}}K_{l_{\dim(\sigma')}}\cdots K_{l_1}\Omega}
{\prod_{\tilde\pi(\rho_k)\subset\sigma}x_k}$$
(where, $\Omega$   is the $d$-form for $\ps$ from \cite[Definition~9.3]{bc})
for the integral generators $e_{l_1},\dots,e_{l_{\dim(\sigma')}}$ 
of a cone $\sigma'\in\Sigma(d-i+\dim\sigma)$ such that
$\tau\subset\sigma'$ and the image of $\sigma'$ is in $\sigma$. This can be seen from the definition of $\Omega$ and 
$\Omega_\sigma$, and
 the fact that the choice of a different basis in the definition of $\Omega$ affects only the sign of the form.
Hence, ${\rm d}f\wedge\Omega_\sigma$ is equal up to a multiple by a constant to
\begin{multline*}
{\rm d}f\wedge\frac{x_{l_1}\cdots x_{l_{\dim(\sigma')}}K_{l_{\dim(\sigma')}}\cdots K_{l_1}\Omega}
{\prod_{\tilde\pi(\rho_k)\subset\sigma}x_k}
\\
=(-1)^{\dim(\sigma')}x_{l_1}\cdots x_{l_{\dim(\sigma')}}
\frac{K_{l_{\dim(\sigma')}}\cdots K_{l_1}(\dd f\wedge\Omega)}
{\prod_{\tilde\pi(\rho_k)\subset\sigma}x_k}
\\
+\sum^{\dim(\sigma')}_{r=1}
(-1)^{r-1} f_{l_r}x_{l_1}\cdots x_{l_{\dim(\sigma')}}\frac{
K_{l_{\dim(\sigma')}}\cdots\widehat{K_{l_r}}\cdots K_{l_1}\Omega}{\prod_{\tilde\pi(\rho_k)\subset\sigma}x_k}.
\end{multline*} 
The first summand is a multiple of $f$ by equation (3) in \cite{m1}, and we claim that
 the second is a multiple of $x_k$ or $\dd x_k$ for  $e_k\in\tau$ such that its image is in the relative 
interior of $\sigma$. Indeed, each summand in $\Omega$ is a multiple of $x_k$ or $\dd x_k$, while the argument in the proof
of Lemma~4.2 in \cite{m2} shows that
$x_{l_r}f_{l_r}$ is divisible by all $x_k$ for  $\tilde\pi(e_k)$ in the relative 
interior of $\sigma$.
Thus, equation (\ref{e:mul}) is proved.

In the second part, 
we need to find a representative in the \v{C}ech cohomology of the image of the form 
$A\Omega_\sigma/f^{p+1}$ by the composition of the maps in (\ref{m:seq}). 
The form $A\Omega_\sigma/f^{p+1}$ has a representation in 
$\check{H}^0({\cal U}|_{V(\tau)},L_{p+1})$ by the \v{C}ech cocycle
$${\fracwithdelims\{\}{A\dd f\wedge K_{i_0}\Omega_\sigma}{f_{i_0}f^{p+1}}}_{i_0},$$
since 
$\dd f\wedge K_{i_0}\Omega_\sigma=f_{i_0}\Omega_\sigma-K_{i_0}(\dd f\wedge\Omega_\sigma)$
and because of equation (\ref{e:mul}).
Note that the maps between cohomologies 
in (\ref{m:seq}) are the connecting homomorphisms of long exact sequences, so that the image of the 
\v{C}ech cocycle representing $A\Omega_\sigma/f^{p+1}$ can be computed
using the following commutative diagrams:
$$
\minCDarrowwidth0.4cm
\begin{CD} 
0 @>>> C^{p+1-s}({\cal U}|_{V(\tau)},L_{s}) @>>> 
C^{p+1-s}({\cal U}|_{V(\tau)},\tilde L_{s}) 
@>\dd>> C^{p+1-s}({\cal U}|_{V(\tau)},L_{s+1}) \\
@.    @AAA                               @AAA             @AAA   \\
0 @>>>C^{p-s}({\cal U}|_{V(\tau)},L_{s}) @>>> 
C^{p-s}({\cal U}|_{V(\tau)},\tilde L_{s}) 
@>\dd>> C^{p-s}({\cal U}|_{V(\tau)},L_{s+1}), 
\end{CD}
$$
where $\tilde L_{s}:=\Omega_{V(\tau)}^{i(\sigma)-p+s-1}(sX_\tau)/\Omega_{V(\tau)}^{i(\sigma)-p+s-1}((s-1)X_\tau)$
and the vertical arrows are the \v{C}ech coboundary maps. When $s=p$,  the cocycle representing $A\Omega_\sigma/f^{p+1}$ 
has a lift to the cochain
$$-\frac{1}{p}{\fracwithdelims\{\}{AK_{i_0}\Omega_\sigma}{f_{i_0}f^{p}}}_{i_0}\in C^{0}({\cal U}|_{V(\tau)},
\Omega_{V(\tau)}^{i(\sigma)-1}(pX_\tau)/\Omega_{V(\tau)}^{i(\sigma)-1}((p-1)X_\tau)),$$
because 
$$
\dd\fracwithdelims()
{AK_{i_0}\Omega_\sigma}{f_{i_0}f^{p}}=
-pA\frac{\dd f\wedge K_{i_0}\Omega_\sigma}{f_{i_0}f^{p+1}}+\frac{\dd(AK_{i_0}\Omega_\sigma)}{f_{i_0}f^{p}}-
\frac{A\dd f_{i_0}\wedge K_{i_0}\Omega_\sigma}{f^{p}f_{i_0}^2}.
$$
The \v{C}ech coboundary of this cochain is
\begin{multline*}
-\frac{1}{p}{\biggl\{\frac{AK_{i_1}\Omega_\sigma}{f_{i_1}f^{p}}-\frac{AK_{i_0}\Omega_\sigma}{f_{i_0}f^{p}}\biggr\}}_{i_0i_1}=
-\frac{1}{p}{\biggl\{\frac{A(f_{i_0}K_{i_1}\Omega_\sigma-f_{i_1}K_{i_0}\Omega_\sigma)}{f_{i_0}f_{i_1}f^{p}}\biggr\}}_{i_0i_1}
\\
=-\frac{1}{p}{\biggl\{\frac{A(-\dd f\wedge K_{i_1}K_{i_0}\Omega_\sigma+K_{i_1}K_{i_0}(\dd f\wedge\Omega_\sigma))}
{f_{i_0}f_{i_1}f^{p}}\biggr\}}_{i_0i_1}
\\
=\frac{1}{p}{\biggl\{\frac{A\dd f\wedge K_{i_1}K_{i_0}\Omega_\sigma}
{f_{i_0}f_{i_1}f^{p}}\biggr\}}_{i_0i_1},
\end{multline*} 
and, by the above diagram, it is a \v{C}ech cocycle in $C^{1}({\cal U}|_{V(\tau)},L_p)$.
Continuing this way, one gets the cocycle 
\begin{multline*}
\frac{(-1)^{p}(-1)^{1+\cdots+p}}{p!}
{\fracwithdelims\{\}{A\dd f\wedge K_{i_{p}}\cdots K_{i_0}\Omega_\sigma}{f_{i_0}\cdots f_{i_{p}}f}}_{i_0\dots i_{p}}
\\
=\frac{(-1)^{i(\sigma)-1+(p(p+1)/2)}}{p!}
{\biggl\{\frac{A K_{i_{p}}\cdots K_{i_0}\Omega_\sigma}{f_{i_0}\cdots f_{i_{p}}}\wedge\frac{\dd f}{f}\biggr\}}_{i_0\dots i_{p}}
\end{multline*} 
in $C^p({\cal U}|_{V(\tau)},\Omega_{V(\tau)}^{i(\sigma)-p}(\log X_\tau))$, which is ready for application of 
the residue map. This produces the desired result.
\end{pf}

To compute the cup product on the residue part of cohomology of semiample hypersurfaces,
we will need to use certain toric residue maps similar to \cite{c2} and \cite{m1}. Let us first
recall some facts about toric residues, and then apply them to our situation.
 
\begin{defn} \cite{c2}
Assume $F_0,\dots,F_d\in S_\beta$ do not vanish simultaneously on a complete $d$-dimensional toric variety $\pp$.
Then the {\it toric residue map} 
$$\res_F: S_\rho/\langle F_0,\dots,F_d\rangle_\rho@>>> {\Bbb C},$$ 
$\rho={(d+1)\beta-\beta_0}$,
is given by the formula $\res_F(H)=\tr_\pp([\varphi_F(H)])$, where 
$\tr_\pp:H^d(\pp,\Omega_{\pp}^d)\rightarrow{\Bbb C}$ is the trace map, and  $[\varphi_F(H)]$  is 
the class represented by the $d$-form ${H\Omega}/({F_0\cdots F_d})$ in the \v{C}ech cohomology with
respect to the open cover $\{x\in\pp:F_j(x)\ne0\}$.
\end{defn}

In \cite[Definition~3.4]{m1} we introduced the following constant.

\begin{defn} \cite{m1}
For $\beta=[\sum_{k=1}^n b_k D_k]\in A_{d-1}(\pp)$ and an ordered  subset $I=\{i_0,\dots,i_d\}$ 
of $\{1,\dots,n\}$, let the constant $c_I^\beta$  be the determinant of the $(d+1)\times(d+1)$ matrix 
obtained from
($\langle m_j, e_{k}\rangle_{1\le j\le d, k\in I})$ by adding the first row 
$(b_{k})_{k\in I}$, where  $m_1,\dots,m_d$ is the same integer basis
of the lattice $M$ as in the definition of $\Omega$. 
\end{defn}

The next result from \cite[Theorem~4.8]{m1} shows how to explicitly compute the toric residue map.

\begin{thm}\label{t:jac} \cite{m1}
Let $\pp$ be a complete toric variety, and let $\beta=[D]\in A_{d-1}(\pp)$ be big and nef.
If $F_0,\dots,F_d\in S_\beta$ do not vanish simultaneously on $\pp$, then:

{\rm(i)} the  toric residue map 
$\res_F: S_\rho/\langle F_0,\dots,F_d\rangle_\rho@>>> {\Bbb C}$,
$\rho={(d+1)\beta-\beta_0}$, is an isomorphism,  

{\rm(ii)} $\res_F(J_F)=d!\vol(\Delta_D)$, where 
$J_F=\det(\partial F_j/\partial x_{i_k})/c_I^\beta \hat{x}_I\in S_{(d+1)\beta-\beta_0}$ 
is the toric Jacobian of  $F_0,\dots,F_d$, and $c_I^\beta\ne0$.
\end{thm}

When $X\subset\pp$ is a nondegenerate hypersurface determined by $f\in S_\beta$, by Lemma~4.10 in \cite{m1},
the polynomials $F_k=x_k\partial f/\partial x_k$, $k\in I$, do not vanish simultaneously on $\pp$
for  $c_I^\beta\ne0$, and 
$J_0(f)=\langle x_k\partial f/\partial x_k:\,k\in I \rangle$.
 Therefore, there is  the corresponding toric residue map
$\res_{F_I}: S_\rho/{J_0(f)}_\rho@>>> {\Bbb C}$. 
However, the diagram (10) in \cite{m1} tells us that this map depends on the choice of $I$, while
$c_I^\beta\res_{F_I}$ doesn't. By the above theorem, $c_I^\beta\res_{F_I}(J)=d!\vol(\Delta_D)$ for
\begin{equation}\label{e:pol}
J=\frac{\det(\partial F_j/\partial x_{k})_{k,j\in I}}{(c_I^\beta)^2 \hat{x}_I}.
\end{equation}

For an $i$-semiample nondegenerate hypersurface $X\subset\ps$ and $\sigma\in\Sigma_X$.
By Proposition~\ref{p:iso}, there is an isomorphism
$$R^\sigma_0(f)_{(i-\dim\sigma+1)\beta-\beta_0+\beta_1^\sigma}\,{\cong}\,
R_0(f_{\sigma'})_{(i-\dim\sigma+1)\beta^{\sigma'}-\beta_0^{\sigma'}}
\,{\cong}\,R_0(f_{\sigma})_{(i-\dim\sigma+1)\beta^{\sigma}-\beta_0^{\sigma}},$$
composition of which we denote by ${\pi_\sigma}_*$.
Since $Y\cap V(\sigma)$ is an ample hypersurface in $V(\sigma)$, determined by $f_{\sigma}$,
we get the corresponding map  
$$c_I^{\beta^\sigma}\res_{{F_\sigma}_I}:R_0(f_{\sigma})_{(i-\dim\sigma+1)\beta^{\sigma}-\beta_0^{\sigma}}\,\cong\,{\Bbb C}.$$ 
It is not difficult to find a polynomial in $R_0(f_{\sigma})_{(i-\dim\sigma+1)\beta^{\sigma}-\beta_0^{\sigma}}$,
which ${\varphi_{\sigma'}}_*$ and ${\pi_\sigma}_*$ map to the polynomials (\ref{e:pol}) 
associated with $f_{\sigma'}\in S(V(\sigma'))_{\beta^{\sigma'}}$ and $f_\sigma\in S(V(\sigma))_{\beta^{\sigma}}$, respectively.
This polynomial is
$$J_{\sigma}=\frac{\det(x_k\partial F_j/\partial x_k)_{k,j\in I}}
{(c_I^{\beta,\sigma})^2\prod_{\tilde\pi(\rho_k)\not\subset\sigma}x_k},$$
where $c_I^{\beta,\sigma}\ne0$ for an  ordered subset $I$ of $\{1,\dots,n\}$
is the determinant of the $(i(\sigma)+1)\times(i(\sigma)+1)$ matrix 
obtained from
$(\langle m^\sigma_j, e_{k}\rangle)_{1\le j\le i(\sigma),k\in I}$ by adding the first row 
$(b_{i})_{i\in I}$ such that $\beta=[\sum_{k=1}^n b_k D_k]$ with $b_k=0$ if $\rho_k\subset\sigma$,
and $m_1^\sigma,\dots,m^\sigma_{i(\sigma)}$ is the same integer basis
of the lattice $M_X\cap\sigma^\perp$ as in the definition of $\Omega_\sigma$. 
One can  verify that  $J_{\sigma}$ does not depend on the choice of $I$, and 
the fact that  $J_{\sigma}$ is indeed a polynomial follows from the arguments at the end of the proof of 
\cite[Proposition~4.1]{c2}.
Using Theorem~\ref{t:jac}, we conclude that
\begin{equation}\label{e:resj}
c_I^{\beta^\sigma}\res_{{F_\sigma}_I}{\pi_\sigma}_*(J_{\sigma})=(i-\dim\sigma)!\,\vol(\Gamma_\sigma),
\end{equation}
where $\Gamma_\sigma$ is the face of $\Delta_D$ corresponding to $\sigma\in\Sigma_X$.

In the case $\dim\sigma=i-1$, we will also need to use a special map constructed the following way.
By equation  (\ref{e:resj}), 
\begin{equation}\label{e:pro}
C\prod_{\tilde\pi(\rho_k)\not\subset\sigma}x_k-\frac{c_I^{\beta^\sigma}\res_{{F_\sigma}_I}{\pi_\sigma}_*
(C\prod_{\tilde\pi(\rho_k)\not\subset\sigma}x_k)}
{\vol(\Gamma_\sigma)}J_{\sigma}\in J^\sigma_0(f)_{2\beta-\beta_0+\beta_1^{\sigma}}
\end{equation}
for $C\in (S/\langle x_k:\,\tilde\pi(\rho_k)\subset\sigma\rangle)_{2\beta-2\beta_0+2\beta_1^{\sigma}}$.
The ideal $J^\sigma_0(f)$ is generated by $J_0(f)$ and $x_k$ for $\tilde\pi(\rho_k)\subset\sigma$,
while the ideal $J_0(f)$ has two minimal generators $f$ and $x_sf_s$ if $\tilde\pi(\rho_s)\not\subset\sigma$.
Define a projection $p_\sigma$ as the composition of the maps:
\begin{multline*}
(S/\langle x_k:\,\tilde\pi(\rho_k)\subset\sigma\rangle)_{2\beta-2\beta_0+2\beta_1^{\sigma}}@>>>
J^\sigma_0(f)_{\rho}@>>> 
J^\sigma_0(f)_{\rho}/\langle f, x_k:\,\tilde\pi(\rho_k)\subset\sigma\rangle)_{\rho}\\
\cong\,
(\langle x_sf_s\rangle/\langle x_k:\,\tilde\pi(\rho_k)\subset\sigma\rangle)_{\rho}
\,\cong\, 
(S/\langle  x_k:\,\tilde\pi(\rho_k)\subset\sigma\rangle)_{\beta-\beta_0+\beta_1^{\sigma}}=
R^\sigma_1(f)_{\beta-\beta_0+\beta_1^{\sigma}},
\end{multline*}
where $\rho=2\beta-\beta_0+\beta_1^{\sigma}$, 
the first map is defined by (\ref{e:pro}) and the fourth map is a division by $x_sf_s/\langle m^\sigma_1, e_{s}\rangle$.
The map $p_\sigma$ is independent of the choice of $s$. Indeed, since
$$\frac{e_{s_1}}{\langle m^\sigma_1, e_{s_1}\rangle}-\frac{e_{s_2}}{\langle m^\sigma_1, e_{s_2}\rangle}$$
belongs to the space spanned by $e_k$ for $\tilde\pi(\rho_k)\subset\sigma$, 
there is an Euler formula (see \cite[Definition~3.9]{bc}) showing that
$$\frac{x_{s_1}f_{s_1}}{\langle m^\sigma_1, e_{s_1}\rangle}-\frac{x_{s_2}f_{s_2}}{\langle m^\sigma_1, e_{s_2}\rangle}$$
belongs to the ideal $\langle f, x_k:\,\tilde\pi(\rho_k)\subset\sigma\rangle$.

For $\sigma\in\Sigma_X$, define  $\res^{\sigma}:R^\sigma_1(f)@>>>{\Bbb C}$  by 0  in all degrees except
$i(\sigma)+1)\beta-2\beta_0+2\beta_1^\sigma$, for which $\res^{\sigma}$ is 
the composition
$$R^\sigma_1(f)_{(i(\sigma)+1)\beta-2\beta_0+2\beta_1^\sigma}@>>>
R^\sigma_0(f)_{(i(\sigma)+1)\beta-\beta_0+\beta_1^\sigma}@>c_I^{\beta^\sigma}\res_{{F_\sigma}_I}{\pi_\sigma}_*
>>{\Bbb C},$$
where the first arrow is a multiplication by 
$\prod_{\tilde\pi(\rho_k)\not\subset\sigma}x_k/((i-\dim\sigma)!\,\vol(\Gamma_\sigma))$.

Now, we are ready to show how to partially split the residue part $H_{\rm res}^s(X\cap V(\tau))$.

\begin{pr}\label{p:pdec}
Let $X\subset\ps$ be an $i$-semiample nondegenerate hypersurface 
defined by $f\in S_\beta$. Given $\tau\in\Sigma$ and  the smallest cone $\sigma\in\Sigma_X$,   containing
the image of $\tau$, then: 
{\rm(i)}
 for $p+q\ne i-\dim\sigma-1$,
$$H_{\rm res}^{p,q}(X\cap V(\tau))=
\sum_{\tau\subset\tau'\in\Sigma(\dim\tau+1)}{\varphi_{\tau,\tau'}}_! H_{\rm res}^{p-1,q-1}(X\cap V(\tau')),
$$ 
 and
$$H_{\rm res}^{p,q}(X\cap V(\tau))\cong R^\sigma_1(f)_{(q+1)\beta-\beta_0+\beta_1^\sigma}
\bigoplus\sum_{\begin{Sb}\tau\subset\tau'\\ \tau'\in\Sigma(\dim\tau+1)\end{Sb}}
{\varphi_{\tau,\tau'}}_! H_{\rm res}^{p-1,q-1}(X\cap V(\tau'))
$$
for $p+q=i-\dim\sigma-1$, where the first factor is included by the residue map,

{\rm(ii)} 
  for  $\dim\sigma<i-1$, $A\in R^\sigma_1(f)_{(p+1)\beta-\beta_0+\beta_1^\sigma}$,
$B\in R^\sigma_1(f)_{(q+1)\beta-\beta_0+\beta_1^\sigma}$,
$$\res_{\tau}(A)\cup\res_{\tau}(B)=
c_p^\sigma\res^\sigma(AB)\varphi_{\tau}^*i^*X^{i(\sigma)-1},$$ 
where $$c_p^\sigma=-(-2\pi\sqrt{-1})^{i(\sigma)-1}\frac{(-1)^{(i(\sigma)-1)(i(\sigma)+2p+2)}}{p!(i(\sigma)-p-1)!},$$
if $p+q=i-\dim\sigma-1$,
and 
$$\res_{\tau}(A)\cup\res_{\tau}(B)=0$$ 
if $p+q\ne i-\dim\sigma-1$,
and, for $\dim\sigma=i-1$,
$$
\res_{\tau}(A)\cup\res_{\tau}(B)=
-\res^\sigma(AB)\varphi_{\tau}^*i^*[\ps]
+\res_{\tau}(p_\sigma(AB)).
$$
\end{pr}

\begin{pf} {\rm(i)}
{From} the exact sequence (\ref{e:exa}) and the vanishing
$$\gr_s^W H^s(X\cap\tta)\cong \gr_s^W H^s((Y\cap\ts)\times({\Bbb C}^*)^{d-i+\dim\sigma-\dim\tau})=0$$
for $s\ne i-\dim\sigma-1$,
by equations (2) and (11) in \cite{m2}, we conclude the first part of the statement.

For the second part, we first use an induction on the dimension of $\tau\in\Sigma$ such that
$\sigma\in\Sigma_X$ is  the smallest cone  containing
the image of $\tau$ to show that 
\begin{equation}\label{e:rvan}
\res_\tau(J^\sigma_1(f)_{(q+1)\beta-\beta_0+\beta_1^\sigma})=0
\end{equation}
in the cohomology of $X\cap V(\tau)$.  
By Proposition~\ref{p:iso} and \cite[Theorem~4.4]{m2}, this is true for $\tau\in\Sigma(d-i+\dim\sigma)$
 such that  $\tilde\pi(\tau)\subset\sigma$, because $X\cap V(\tau)$ in $V(\tau)$ is big and nef.
By induction, assume that (\ref{e:rvan}) holds for $\tau'\in\Sigma(k)$ such that
$\sigma\in\Sigma_X$ is  the smallest cone  containing $\tilde\pi(\tau')$. We need to prove (\ref{e:rvan}) for
$\tau\in\Sigma(k-1)$ satisfying the same condition. 
Let $A\in J^\sigma_1(f)_{(q+1)\beta-\beta_0+\beta_1^\sigma}$, then, by Poincar\'e duality, 
$\res_\tau(A)=0$  in $H^{p,q}(X\cap V(\tau))$, $p=i-\dim\sigma-1-q$, if and only if 
$\int_X$ of the cup product of $\res_\tau(A)$ with all elements in 
$H^{d(\tau)-1-p,d(\tau)-1-q}(X\cap V(\tau))$ vanishes
(note that $X$ is not necessarily connected).
{From} Section~\ref{s:ar} we know that  
$$H^*(X\cap V(\tau))=H^*_{\rm toric}(X\cap V(\tau))\oplus H^*_{\rm res}(X\cap V(\tau)),$$
and 
\begin{equation}\label{e:vcup}
\int_X H^*_{\rm toric}(X\cap V(\tau))\cup H^*_{\rm res}(X\cap V(\tau))=0.
\end{equation}
On the other hand, for $\dim\tau\ne d-i+\dim\sigma$,
$$H_{\rm res}^{d(\tau)-1-p,d(\tau)-1-q}(X\cap V(\tau))=
\sum_{\tau\subset\tau'\in\Sigma(\dim\tau+1)}{\varphi_{\tau,\tau'}}_! 
H_{\rm res}^{d(\tau)-2-p,d(\tau)-2-q}(X\cap V(\tau'))
$$ 
by the first part of this proposition. Since
$$\res_\tau(A)\cup{\varphi_{\tau,\tau'}}_!h={\varphi_{\tau,\tau'}}_!(\varphi_{\tau,\tau'}^*\res_\tau(A)\cup h)$$ 
by the Gysin projection formula, it is enough to show  $\varphi_{\tau,\tau'}^*(\res_\tau(A))=0$ for
$\tau\subset\tau'\in\Sigma(\dim\tau+1)$.
 If $\tilde\pi(\tau')\not\subset\sigma$, then there is $\rho_l\subset\tau'$ such that  
$\tilde\pi(\rho_l)\not\subset\sigma$. In this case, by Proposition~\ref{p:coc},
$\varphi_{\tau,\tau'}^*\res_\tau(A)$ is represented by the restriction 
of the \v{C}ech cocycle
$$\frac{(-1)^{i(\sigma)-1+(q(q+1)/2)}}{q!}
\fracwithdelims\{\}{AK_{i_{q}}\cdots K_{i_0}\Omega_\sigma}{f_{i_0}\cdots f_{i_{q}}}_{i_0\dots i_{q}}
$$ 
to $X\cap V(\tau')$. But each term in the form $K_{i_{q}}\cdots K_{i_0}\Omega_\sigma$ on a nonempty open set 
$U_{i_{q}}\cap\cdots\cap U_{i_0}\cap V(\tau')$ contains $\dd x_l$ or $x_l$ and, 
therefore, vanishes on $X\cap V(\tau')$.  
If $\tilde\pi(\tau')\subset\sigma$, then, by  the induction assumption,
$\res_{\tau'}(A)=0$   in $H^{p,q}(X\cap V(\tau'))$. The commutative diagram
$$ 
\begin{CD}
S_{(q+1)\beta-\beta_0+\beta_1^\sigma}@>\res_\tau>>H^q(X\cap V(\tau),\Omega^{i(\sigma)-1-q}_{X\cap V(\tau)})\\
@|  @VV{\scriptstyle{\scriptstyle\varphi_{\tau,\tau'}^*}}V  \\
S_{(q+1)\beta-\beta_0+\beta_1^\sigma}@>\res_{\tau'}>>H^q(X\cap V(\tau'),\Omega^{i(\sigma)-1-q}_{X\cap V(\tau')}).
\end{CD}
$$ 
gives  
$\varphi_{\tau,\tau'}^*\res_\tau(A)=0$ again. This finishes the proof of  (\ref{e:rvan}) by induction, which
shows that the map
$$R^\sigma_1(f)_{(q+1)\beta-\beta_0+\beta_1^\sigma}@>\res_\tau>>H_{\rm res}^{p,q}(X\cap V(\tau))$$
is well defined.

Finally, for $\sigma'\in\Sigma(d-i+\dim\sigma)$
 such that  $\tilde\pi(\sigma')\subset\sigma$, we have a  diagram
\begin{equation}\label{e:che}
\minCDarrowwidth0.5cm
\begin{CD}
R^\sigma_1(f)_{(q+1)\beta-\beta_0+\beta_1^\sigma}@>\res_\tau>>H_{\rm res}^{p,q}(X\cap V(\tau))@>>>
H^{p,q}\bigl(PH^{p+q}(X\cap\tta)\bigr) \\
@VVV @VVV @.  \\
R_1(f_{\sigma'})_{(q+1)\beta^{\sigma'}-\beta_0^{\sigma'}}@>\res_{\sigma'}>>
H_{\rm res}^{p,q}(X\cap V({\sigma'}))@>>>H^{p,q}\bigl(PH^{p+q}(X\cap\ttt_{\sigma'})\bigr)
\end{CD}\end{equation}
and, by equation (12) in \cite{m2}, the isomorphisms
$$H^{p,q}\bigl(PH^{p+q}(X\cap\tta)\bigr)\,\cong\, H^{p,q}\bigl(PH^{p+q}(\pi(X)\cap\ts)\bigr)
\,\cong\, H^{p,q}\bigl(PH^{p+q}(X\cap\ttt_{\sigma'})\bigr)$$
induced by $\pi:\ps@>>>\psx$. Since $X\cap V({\sigma'})$ in $V({\sigma'})$ is big and nef,
the composition on the bottom of $(\ref{e:che})$ is an isomorphism, by Theorem~4.4 in \cite{m1}.
Hence, the same holds for the composition on the top of $(\ref{e:che})$. Now, the exact sequence  (\ref{e:exa})
produces the result.

{\rm(ii)} Let $\dim\sigma<i-1$. As in part {\rm(i)}, we prove the formula
\begin{equation}\label{e:cupi}
\res_{\tau}(A)\cup\res_{\tau}(B)=
c_p^\sigma\res^\sigma(AB)\varphi_{\tau}^*i^*X^{i(\sigma)-1}
\end{equation}
by an induction on the dimension of $\tau\in\Sigma$ such that
$\sigma\in\Sigma_X$ is  the smallest cone  containing
the image of $\tau$.
When $\dim\tau=d-i+\dim\sigma$, 
the hypersurface $X\cap V(\tau)$ in $V(\tau)$  is big and nef. Therefore, by 
the results in \cite[page~104]{m1} and the definition of $\res^\sigma$, 
$$\int_{X\cap V(\tau)}\res_{\tau}(A)\cup\res_{\tau}(B)=c_p^\sigma(i-\dim\sigma)!\vol(\Gamma_\sigma)\res^\sigma(AB)$$
(the  formula we use in \cite[page~104]{m1} must be corrected by
a factor of $2\pi\sqrt{-1}$).
By \cite[Section~5.3]{f1}, $(X\cap V(\tau))^{i(\sigma)}=(i-\dim\sigma)!\vol(\Gamma_\sigma)$,
whence  (\ref{e:cupi}) holds in this case. Assume that (\ref{e:cupi}) holds for $\dim\tau=k$, 
we will prove (\ref{e:cupi})
for $\dim\tau=k-1$. By Poincar\'e duality, 
it suffices to show that the cup product of both sides of (\ref{e:cupi}) with an element of 
$H^{d(\tau)-i(\sigma),d(\tau)-i(\sigma)}(X\cap V(\tau))$
 produces the same result. The cohomology of $X\cap V(\tau)$ decomposes into the toric and residue parts.
Using the relations in $H^*_{\rm toric}(X\cap V(\tau))$, coming from the cohomology of the ambient toric variety,
we can assume that an element from the toric part is $\varphi_{\tau}^*i^*[V(\gamma)]$ with $\gamma\in\Sigma(\dim\sigma-\dim\tau)$
such that $\gamma\cap\tau=0$. Then we have 
$$\varphi_{\tau}^*i^*[V(\gamma)]=
\alpha_{\gamma\tau}[X\cap V(\tau+\gamma)],
$$
where
$\alpha_{\gamma\tau}=\frac{\mult(\gamma)\mult(\tau)}{\mult(\tau+\gamma)}$
in $H^*(X\cap V(\tau))$, 
because 
$$[V(\tau)]\cdot[V(\gamma)]=\alpha_{\gamma\tau}[X\cap V(\tau+\gamma)]$$
in the Chow ring $A^*(\ps)$. 
Hence, 
$$
\varphi_{\tau}^*i^*[V(\gamma)]\cup\res_{\tau}(A)\cup\res_{\tau}(B)=
\alpha_{\gamma\tau}{\varphi_{\tau+\gamma,\tau}}_!(\varphi_{\tau+\gamma,\tau}^*\res_{\tau}(A)\cup
\varphi_{\tau+\gamma,\tau}^*\res_{\tau}(B))
$$
and 
\begin{multline*}\varphi_{\tau}^*i^*[V(\gamma)]\cup\varphi_{\tau}^*i^*X^{i(\sigma)-1}=
\alpha_{\gamma\tau}{\varphi_{\tau+\gamma,\tau}}_!\varphi_{\tau+\gamma,\tau}^*\varphi_{\tau}^*i^*X^{i(\sigma)-1}
\\
=\alpha_{\gamma\tau}{\varphi_{\tau+\gamma,\tau}}_!\varphi_{\tau+\gamma}^*i^*X^{i(\sigma)-1}.
\end{multline*}
We can  assume that $\tau+\gamma$ forms a cone in $\Sigma$,
and 
that $\gamma$ lies in $\sigma$, because, otherwise, $\varphi_{\tau+\gamma,\tau}^*\res_{\tau}(A)=0$
by the same argument as in \cite[page~102]{m1}, and 
${\varphi_{\tau+\gamma,\tau}}_!\varphi_{\tau+\gamma}^*i^*X^{i(\sigma)-1}=0$ by Lemma~\ref{l:int}.
In this case,
\begin{multline*}
\varphi_{\tau+\gamma,\tau}^*\res_{\tau}(A)\cup
\varphi_{\tau+\gamma,\tau}^*\res_{\tau}(B)=\res_{\tau+\gamma}(A)\cup\res_{\tau+\gamma}(B)
\\=
c_p^\sigma\res^\sigma(AB)\varphi_{\tau+\gamma}^*i^*X^{i(\sigma)-1}
\end{multline*}
by the induction assumption.
To finish the induction proof of (\ref{e:cupi}), we also need to consider the elements 
$h\in H_{\rm res}^*(X\cap V(\tau))$, the cup product of which vanishes with the right hand side
of  (\ref{e:cupi}), by equation (\ref{e:vcup}). The same should hold for the left hand side of  (\ref{e:cupi}).
By part 
{\rm(i)} of this proposition, there are two possibilities: 
$h={\varphi_{\tau,\tau'}}_!h'$ with $h'\in H_{\rm res}^*(X\cap V(\tau'))$ or $h=\res_\tau(C)$.
In the first case,
$${\varphi_{\tau,\tau'}}_!h'\cup\res_{\tau}(A)\cup\res_{\tau}(B)={\varphi_{\tau,\tau'}}_!(h'\cup
\varphi^*_{\tau,\tau'}\res_{\tau}(A)\cup\varphi^*_{\tau,\tau'}\res_{\tau}(B)).$$
If  $\tilde\pi(\tau')\not\subset\sigma$, then $\varphi^*_{\tau,\tau'}\res_{\tau}(A)=0$ as above.
If  $\tilde\pi(\tau')\subset\sigma$, then we can apply the induction assumption to deduce that
$$\varphi^*_{\tau,\tau'}\res_{\tau}(A)\cup\varphi^*_{\tau,\tau'}\res_{\tau}(B)
=\res_{\tau'}(A)\cup\res_{\tau'}(B)$$
belongs to $H_{\rm toric}^*(X\cap V(\tau'))$, the cup product of which with $h'$ is zero.
In the second case, we have 
$$\res_{\tau}(A)\cup\res_{\tau}(B)\cup\res_\tau(C)\in H^{d(\tau)-1,d(\tau)-1}(X\cap V(\tau))$$
where $d(\tau)=d-\dim\tau$. Up to a constant multiple, the cocycle that represents this cup product is
$$\biggl\{ABC\frac{K_{i_p}\cdots K_{i_0}\Omega}{f_{i_0}\cdots f_{i_p}}\wedge
\frac{K_{i_{i(\sigma)-1}}\cdots K_{i_p}\Omega}{f_{i_p}\cdots f_{i_{i(\sigma)-1}}}
\wedge
\frac{K_{i_{d(\tau)-1}}\cdots K_{i_{i(\sigma)-1}}\Omega}{f_{i_{i(\sigma)-1}}\cdots f_{i_{d(\tau)-1}}}\biggr\}
_{i_0\dots i_{d(\tau)-1}}$$
in $C^{d(\tau)-1}({\cal U}|_{X\cap V(\tau)},\Omega_{X\cap V(\tau)}^{d(\tau)-1})$.
It follows from \cite[Proposition~5.3]{c2} that
there are $i-\dim\sigma$ open sets from the open cover ${\cal U}|_{X\cap V(\tau)}$ which cover
$X\cap V(\tau)$. Passing to this refinement, we conclude that the above cup product vanishes, since
$d-\dim\tau>i-\dim\sigma$. Thus, the equation (\ref{e:cupi}) is shown.
The proof  of $\res_{\tau}(A)\cup\res_{\tau}(B)=0$,
for $p+q\ne i-\dim\sigma-1$, follows the same steps as above.
 
When $\dim\sigma=i-1$, we will use a different strategy.
The cup product $\res_{\tau}(A)\cup\res_{\tau}(B)$ is represented by the \v Cech cocycle
$$
{\fracwithdelims\{\}{AB\langle m_1^\sigma,e_{i_0}\rangle^2(\prod_{\tilde\pi(\rho_k)\not\subset\sigma}x_k)^2}
{(x_{i_0}f_{i_0})^2}}_{i_0}.
$$
Write 
$$AB\prod_{\tilde\pi(\rho_k)\not\subset\sigma}x_k=P+\res^\sigma(AB)J_{\sigma},$$
where $P$ is the polynomial determined by (\ref{e:pro}) with $C=AB$.
Using the definition of the projection 
$$p_\sigma:(S/\langle x_k:\,\tilde\pi(\rho_k)\subset\sigma\rangle)_{2\beta-2\beta_0+2\beta_1^{\sigma}}@>>>
R^\sigma_1(f)_{\beta-\beta_0+\beta_1^{\sigma}}$$
before this proposition, one can see that  the cocycle
$$
{\fracwithdelims\{\}{P\langle m_1^\sigma,e_{i_0}\rangle^2\prod_{\tilde\pi(\rho_k)\not\subset\sigma}x_k}
{(x_{i_0}f_{i_0})^2}}_{i_0}.
$$
is equal to
$$
{\fracwithdelims\{\}{p_\sigma(AB)\langle m_1^\sigma,e_{i_0}\rangle\prod_{\tilde\pi(\rho_k)\not\subset\sigma}x_k}
{x_{i_0}f_{i_0}}}_{i_0},
$$
because the difference is by multiples of $f$ and $x_k$, such that $\tilde\pi(\inte\rho_k)\subset\inte\sigma$ and
$\rho_k\subset\tau$, which vanish on $X\cap V(\tau)$.
By Proposition~\ref{p:coc}, this cocycle represents $\res_\tau(p_\sigma(AB))$.
Therefore,  $\res_{\tau}(A)\cup\res_{\tau}(B)-\res_\tau(p_\sigma(AB))$ is represented by the cocycle
\begin{equation}\label{e:cocd}
{\fracwithdelims\{\}{\res^\sigma(AB)J_{\sigma}\langle m_1^\sigma,e_{i_0}\rangle^2\prod_{\tilde\pi(\rho_k)\not\subset\sigma}x_k}
{(x_{i_0}f_{i_0})^2}}_{i_0}.
\end{equation}
As above, we can change $J_{\sigma}$ by multiples of $f$ and $x_k$, such that $\tilde\pi(\inte\rho_k)\subset\inte\sigma$ and
$\rho_k\subset\tau$, because they vanish on $X\cap V(\tau)$. Also, we can assume that the index $i_0$ takes values in
$\{s_1,s_2\}$ for $c_{s_1s_2}^{\beta,\sigma}\ne 0$, since $x_{s_1}f_{s_1}$ and $x_{s_2}f_{s_2}$ do not
vanish simultaneously on $V(\tau)$, by  \cite[Lemma~4.10]{m1}.
Using an Euler formula, similar to the one in the proof of Proposition~4.7 in \cite{m1},
we get that 
$$c_{s_1s_2}^{\beta,\sigma}f-\langle m_1^\sigma,e_{s_2}\rangle x_{s_1}f_{s_1}+
\langle m_1^\sigma,e_{s_1}\rangle x_{s_2}f_{s_2}$$
with $c_{s_1s_2}^{\beta,\sigma}\ne 0$ 
is a multiple of  $x_k$, such that $\tilde\pi(\inte\rho_k)\subset\inte\sigma$ and
$\rho_k\subset\tau$. Hence, modulo these multiples and multiples of $f$,
\begin{align*}
J_{\sigma}&=\frac{\det\begin{pmatrix}
(x_{s_1}\partial_{s_1})^2f& (x_{s_1}\partial_{s_1})(x_{s_2}\partial_{s_2})f\\
(x_{s_2}\partial_{s_2})(x_{s_1}\partial_{s_1})f& (x_{s_2}\partial_{s_2})^2f
\end{pmatrix}}
{(c_{s_1s_2}^{\beta,\sigma})^2\prod_{\tilde\pi(\rho_k)\not\subset\sigma}x_k}
\\
&\equiv
\frac{\det\begin{pmatrix}
x_{s_1}f_{s_1}& x_{s_2}f_{s_2}\\
x_{s_1}x_{s_2}f_{s_1s_2}& x_{s_2}f_{s_2}+x^2_{s_2}f_{s_2s_2}
\end{pmatrix}}
{\langle m_1^\sigma,e_{s_2}\rangle c_{s_1s_2}^{\beta,\sigma}\prod_{\tilde\pi(\rho_k)\not\subset\sigma}x_k}
\equiv
\frac{\det\begin{pmatrix}
f & x_{s_2}f_{s_2}\\
x_{s_2}f_{s_2}& x_{s_2}f_{s_2}+x^2_{s_2}f_{s_2s_2}
\end{pmatrix}}
{\langle m_1^\sigma,e_{s_2}\rangle^2\prod_{\tilde\pi(\rho_k)\not\subset\sigma}x_k}
\\
&\equiv
-\frac{(x_{s_2}f_{s_2})^2}{\langle m_1^\sigma,e_{s_2}\rangle^2\prod_{\tilde\pi(\rho_k)\not\subset\sigma}x_k}
\equiv
-\frac{(x_{s_1}f_{s_1})^2}{\langle m_1^\sigma,e_{s_1}\rangle^2\prod_{\tilde\pi(\rho_k)\not\subset\sigma}x_k}.
\end{align*}
Therefore, the cocycle (\ref{e:cocd}) coincides with 
$
\{-\res^\sigma(AB)\}_{i_0},
$
which comes from the restriction
of the global constant function $-\res^\sigma(AB)\in H^0(V(\tau),{\cal O}_{V(\tau)})\cong{\Bbb C}$.
We have a commutative diagram:
$$
\minCDarrowwidth0.6cm
\begin{CD}
H^0(X\cap V(\tau))@.\cong@. H^0(X\cap V(\tau),{\cal O}_{X\cap V(\tau)})\\
@Ai^*AA @. @Ai^*AA  \\
H^0(V(\tau)) @.\cong@. H^0(V(\tau),{\cal O}_{V(\tau)}).
\end{CD}
$$ 
The isomorphism on the bottom of this diagram sends the fundamental class $[V(\tau)]$ to $1$,
whence the isomorphism  on the top sends $-\res^\sigma(AB)i^*[V(\tau)]$ to the element represented by the 
cocycle (\ref{e:cocd}). This finishes the proof of the proposition.
\end{pf}

After all the intermediate technical steps we can 
describe the generators which span the residue part of the cohomology of semiample nondegenerate hypersurfaces.

\begin{pr}\label{p:pdes} Let $X\subset\ps$ be an $i$-semiample nondegenerate hypersurface 
defined by $f\in S_\beta$. Then
$$H_{\rm res}^{p,q}(X)\cong\sum_{\begin{Sb}\sigma\in\Sigma_X \\ 
\gamma\in\Sigma((\dim\sigma+p+q+1-i)/2)\\ 
\tilde\pi(\inte\gamma)\subset\inte\sigma\end{Sb}}{\varphi_{\gamma}}_!\res_\gamma 
R^\sigma_1(f)_{(q+1-\dim\gamma)\beta-\beta_0+\beta_1^\sigma}.$$
\end{pr}

\begin{pf} We can actually prove by induction  on the codimension of $\tau\in\Sigma$ a more general statement: 
$$H_{\rm res}^{p,q}(X\cap V(\tau))\cong\sum_{\begin{Sb} \tilde\pi(\tau)\subset\sigma\in\Sigma_X \\
\tau\subset\gamma\in\Sigma((p+q+1-i(\sigma))/2+\dim\tau)\\ 
\tilde\pi(\inte\gamma)\subset\inte\sigma\end{Sb}}\hspace{-1.2cm}{\varphi_{\tau,\gamma}}_!\res_\gamma 
R^\sigma_1(f)_{(q+1-\dim\gamma+\dim\tau)\beta-\beta_0+\beta_1^\sigma}.$$
This is certainly true for the maximally dimensional $\tau\in\Sigma(d)$, because $X\cap V(\tau)$ is empty by Proposition~1.6 and
Remark~2.2 in \cite{m2} and there is no $\gamma$ satisfying the conditions on the right side of the formula.
Assume that the formula holds for $\tau\in\Sigma(k-1)$.  
Then Proposition~\ref{p:pdec} and the induction assumption  produce the result for $\tau\in\Sigma(k)$. 
\end{pf}

The next result describes the product structure 
on the residue part of  the cohomology of semiample nondegenerate hypersurfaces.

\begin{pr}\label{p:produ}
Let $X\subset\ps$ be an $i$-semiample nondegenerate hypersurface 
defined by $f\in S_\beta$,  then:

{\rm(i)}  for $\sigma\in\Sigma_X$ with $\dim\sigma<i-1$, $\gamma_1,\gamma_2\in\Sigma$,
such that $\tilde\pi(\inte\gamma_k)\subset\inte\sigma$ with $k=1,2$, and $A\in R^\sigma_1(f)_{(p+1)\beta-\beta_0+\beta_1^\sigma}$,
$B\in R^\sigma_1(f)_{*\beta-\beta_0+\beta_1^\sigma}$,  
$$\frac{{\varphi_{\gamma_1}}_!\res_{\gamma_1}(A)}{\mult(\gamma_1)}\cup\frac{{\varphi_{\gamma_2}}_!\res_{\gamma_2}(B)}
{\mult(\gamma_2)}=
c_p^\sigma\res^\sigma(AB)i^*X^{i(\sigma)-1}\prod_{\rho_k\subset\gamma_1}i^*D_k\prod_{\rho_k\subset\gamma_2}i^*D_k,$$
where $i(\sigma)=i-\dim\sigma$ and 
$$c_p^\sigma=-(-2\pi\sqrt{-1})^{i(\sigma)-1}\frac{(-1)^{(i(\sigma)-1)(i(\sigma)+2p+2)}}{p!(i(\sigma)-p-1)!},$$

{\rm(ii)}  for $\dim\sigma=i-1$,
$A,B$ and $\gamma_1,\gamma_2$, as in part {\rm(i)},
\begin{multline*}
\frac{{\varphi_{\gamma_1}}_!\res_{\gamma_1}(A)}{\mult(\gamma_1)}\cup\frac{{\varphi_{\gamma_2}}_!\res_{\gamma_2}(B)}
{\mult(\gamma_2)}=
-\res^\sigma(AB)i^*[\ps]\prod_{\rho_k\subset\gamma_1}i^*D_k\prod_{\rho_k\subset\gamma_2}i^*D_k\\
+\prod_{\rho_k\subset\gamma_1\cap\gamma_2}i^*D_k\frac{{\varphi_{\gamma_1+\gamma_2}}_!\res_{\gamma_1+\gamma_2}(p_\sigma(AB))}
{\mult(\gamma_1+\gamma_2)},
\end{multline*}

{\rm(iii)} for  $\gamma_1,\gamma_2\in\Sigma$  satisfying neither part {\rm(i)} nor {\rm(ii)},
$${\varphi_{\gamma_1}}_!\res_{\gamma_1}(A)\cup{\varphi_{\gamma_2}}_!\res_{\gamma_2}(B)=0.$$
\end{pr}

\begin{pf}
 Applying the Gysin projection formula, we get
\begin{align*}
{\varphi_{\gamma_1}}_!&\res_{\gamma_1}(A)\cup{\varphi_{\gamma_2}}_!\res_{\gamma_2}(B)
\\
&=
{\varphi_{\gamma_1\cap\gamma_2}}_!{\varphi_{\gamma_1\cap\gamma_2,\gamma_1}}_!
\res_{\gamma_1}(A)\cup{\varphi_{\gamma_1\cap\gamma_2}}_!{\varphi_{\gamma_1\cap\gamma_2,\gamma_2}}_!\res_{\gamma_2}(B)
\\
&={\varphi_{\gamma_1\cap\gamma_2}}_!(\varphi_{\gamma_1\cap\gamma_2}^*{\varphi_{\gamma_1\cap\gamma_2}}_!
{\varphi_{\gamma_1\cap\gamma_2,\gamma_1}}_!
\res_{\gamma_1}(A)\cup{\varphi_{\gamma_1\cap\gamma_2,\gamma_2}}_!\res_{\gamma_2}(B))
\\
&=
{\varphi_{\gamma_1\cap\gamma_2}}_!\varphi_{\gamma_1\cap\gamma_2}^*{\varphi_{\gamma_1\cap\gamma_2}}_!
({\varphi_{\gamma_1\cap\gamma_2,\gamma_1}}_!\res_{\gamma_1}(A)\cup
{\varphi_{\gamma_1\cap\gamma_2,\gamma_2}}_!\res_{\gamma_2}(B))
\\
&={\varphi_{\gamma_1\cap\gamma_2}}_!\varphi_{\gamma_1\cap\gamma_2}^*{\varphi_{\gamma_1\cap\gamma_2}}_!
{\varphi_{\gamma_1\cap\gamma_2,\gamma_2}}_!
(\varphi_{\gamma_1\cap\gamma_2,\gamma_2}^*
{\varphi_{\gamma_1\cap\gamma_2,\gamma_1}}_!\res_{\gamma_1}(A)\cup
\res_{\gamma_2}(B))
\\
&={\varphi_{\gamma_1\cap\gamma_2}}_!\varphi_{\gamma_1\cap\gamma_2}^*{\varphi_{\gamma_2}}_!
(\alpha_{\gamma_1\gamma_2}{\varphi_{\gamma_2,\gamma_1+\gamma_2}}_!\varphi_{\gamma_1,\gamma_1+\gamma_2}^*
\res_{\gamma_1}(A)\cup
\res_{\gamma_2}(B))
\\
&=\alpha_{\gamma_1\gamma_2}
{\varphi_{\gamma_1\cap\gamma_2}}_!\varphi_{\gamma_1\cap\gamma_2}^*{\varphi_{\gamma_2}}_!
{\varphi_{\gamma_2,\gamma_1+\gamma_2}}_!(\varphi_{\gamma_1,\gamma_1+\gamma_2}^*
\res_{\gamma_1}(A)\cup\varphi_{\gamma_2,\gamma_1+\gamma_2}^*\res_{\gamma_2}(B))
\\
&=\alpha_{\gamma_1\gamma_2}
{\varphi_{\gamma_1\cap\gamma_2}}_!\varphi_{\gamma_1\cap\gamma_2}^*
{\varphi_{\gamma_1+\gamma_2}}_!(\varphi_{\gamma_1,\gamma_1+\gamma_2}^*
\res_{\gamma_1}(A)\cup\varphi_{\gamma_2,\gamma_1+\gamma_2}^*\res_{\gamma_2}(B))
\\
&=\alpha_{\gamma_1\gamma_2}i^*[V(\gamma_1\cap\gamma_2)]
{\varphi_{\gamma_1+\gamma_2}}_!(\varphi_{\gamma_1,\gamma_1+\gamma_2}^*
\res_{\gamma_1}(A)\cup\varphi_{\gamma_2,\gamma_1+\gamma_2}^*\res_{\gamma_2}(B)),
\end{align*}
where we also used the commutative diagram (see \cite[Lemma~5.4]{m1}):
$$
\minCDarrowwidth0.6cm
\begin{CD}
H^*(X\cap V(\gamma_1))@>{\varphi_{\gamma_1\cap\gamma_2,\gamma_1}}_!>>H^*(X\cap V(\gamma_1\cap\gamma_2))\\
@V\varphi_{\gamma_1,\gamma_1+\gamma_2}^*VV @V\varphi_{\gamma_1\cap\gamma_2}^*VV \\
H^*(X\cap V(\gamma_1+\gamma_2))@>\alpha_{\gamma_1\gamma_2}{\varphi_{\gamma_2,\gamma_1+\gamma_2}}_!>> H^*(X\cap V(\gamma_2)), 
\end{CD}
$$
with $\alpha_{\gamma_1\gamma_2}$ satisfying the equality 
$[X\cap V(\gamma_1)]\cup[X\cap V(\gamma_2)]=\alpha_{\gamma_1\gamma_2}[X\cap V(\gamma_1+\gamma_2)]$ in 
$H^*(X\cap V(\gamma_1\cap\gamma_2))$.

Note that the cup product calculated  above is zero, if $\gamma_1$ and $\gamma_2$ do not span a cone in $\Sigma$.
It also vanishes when $\gamma_1$ and $\gamma_2$ have different smallest cones $\sigma_1$ and $\sigma_2$ in $\Sigma_X$
containing them. Indeed, we can assume that $\sigma_1$ is not contained in $\sigma_2$.  Then
$\varphi_{\gamma_1,\gamma_1+\gamma_2}^*\res_{\gamma_1}(A)$  vanishes by the same argument as in \cite[p.~102]{m1}.

To determine the constant $\alpha_{\gamma_1\gamma_2}$ first notice that
$[X\cap V(\gamma)]=i^*[V(\gamma)]$, for $\gamma\in\Sigma$, by the arguments in the proof of Lemma~5.7 in \cite{m1}.
On the other hand, we can write $\gamma_k=\gamma_k'+\gamma_1\cap\gamma_2$, for $k=1,2$, so that $\gamma'\in\Sigma$ and
$\dim\gamma_k=\dim\gamma_k'+\dim(\gamma_1\cap\gamma_2)$.
In this case, 
$$[V(\gamma)]=\frac{\mult(\gamma)}{\mult(\gamma')\mult(\gamma_1\cap\gamma_2)}
\varphi_{\gamma_1\cap\gamma_2}^*[V(\gamma')]$$
in $H^*(V(\gamma_1\cap\gamma_2))$, where $\gamma$ is $\gamma_1$, $\gamma_2$ or $\gamma_1+\gamma_2$ and 
$\gamma'$ is $\gamma_1'$, $\gamma_2'$ or $\gamma_1'+\gamma_2'$, respectively.
This is  because
$$V(\gamma)=\frac{\mult(\gamma)}{\mult(\gamma')\mult(\gamma_1\cap\gamma_2)}V(\gamma_1\cap\gamma_2)\cdot V(\gamma')$$
in the Chow ring $A^*(\ps)$ for the corresponding $\gamma$  and $\gamma'$. 
Since       
$$V(\gamma_1')\cdot V(\gamma_2')=\frac{\mult(\gamma_1')\mult(\gamma_2')}{\mult(\gamma_1'+\gamma_2')}
V(\gamma_1'+\gamma_2'),$$
it is not difficult to calculate
$$\alpha_{\gamma_1\gamma_2}=\frac{\mult(\gamma_1)\mult(\gamma_2)}
{\mult(\gamma_1\cap\gamma_2)\mult(\gamma_1+\gamma_2)}.$$
Hence, 
\begin{multline*}
{\varphi_{\gamma_1}}_!\res_{\gamma_1}(A)\cup{\varphi_{\gamma_2}}_!\res_{\gamma_2}(B)
\\=
\frac{\mult(\gamma_1)\mult(\gamma_2)}{\mult(\gamma_1+\gamma_2)}
\prod_{\rho_k\subset\gamma_1\cap\gamma_2}i^*D_k
{\varphi_{\gamma_1+\gamma_2}}_!(\res_{\gamma_1+\gamma_2}(A)\cup\res_{\gamma_1+\gamma_2}(B))
\end{multline*}
Applying Proposition~\ref{p:pdes}, we get the result.
\end{pf}

To give a very explicit description of the cohomology of semiample  hypersurfaces we introduce the following rings.

\begin{defn} 
Given a semiample class $[D]\in A_{d-1}(\ps)$   and $\sigma\in\Sigma_D$, 
let 
$$U^\sigma(D)=\biggl\langle \prod_{\rho_k\subset\gamma\in\Sigma}D_k: \tilde\pi(\inte\gamma)\subset\inte\sigma\biggr\rangle$$
be the ideal in ${\Bbb C}[D_1,\dots,D_n]$.
Define the bigraded ring  
$$A_1^\sigma(D)_{*,*}=U^\sigma(D)/\{u\in U^\sigma(D):\,uvX^{i(\sigma)}\in(P(\Sigma)+SR(\Sigma))\text{ for all }v\in U^\sigma(D)\},$$ 
where $D_k$ has the degree (1,1).
\end{defn} 

Here is our main result. 

\begin{thm}\label{t:main} Let $X\subset\ps$ be an $i$-semiample nondegenerate hypersurface 
defined by $f\in S_\beta$.
Then there is a ring isomorphism
$$\bigoplus_{p,q}H^{p,q}(X)\cong\bigoplus_{p,q} A_1(X)_{p,q}\oplus\biggl(\bigoplus_{\sigma\in\Sigma_X}
A_1^\sigma(X)_{s,s}
\otimes R^\sigma_1(f)_{(q-s+1)\beta-\beta_0+\beta_1^\sigma}\biggr),$$
where $s=(p+q-i-\dim\sigma+1)/2$.  
The product structure on the right side is given by:\\
{\rm(a)}   $a\cdot b=ab$ for $a,b\in A_1(X)$,\\
{\rm(b)} $a\cdot (u\otimes g)= (au)\otimes g${ for } $a\in A_1(X)$, $u\in A_1^\sigma(X)$, 
$g\in R^\sigma_1(f)_{*\beta-\beta_0+\beta_1^\sigma}$,\\
{\rm(c)} $(u\otimes g)\cdot (v\otimes h)=c_r^\sigma \res^\sigma(g h) X^{i(\sigma)-1}uv$
for $\dim\sigma<i-1$, $u, v\in A_1^\sigma(X)$, 
$g\in R^\sigma_1(f)_{(r+1)\beta-\beta_0+\beta_1^\sigma}$, $h\in R^\sigma_1(f)_{t\beta-\beta_0+\beta_1^\sigma}$
($\res^\sigma(g h)$ vanishes unless $r+t=i(\sigma)$,  $i(\sigma)=i-\dim\sigma$),\\
{\rm(d)} $(u\otimes g)\cdot (v\otimes h)=- \res^\sigma(g h) uv+(uv)\otimes p_\sigma(gh)$
for $\dim\sigma=i-1$, $u$, $v$, $g$, $h$, as in part {\rm(b)},\\
{\rm(e)} $(u\otimes g)\cdot (v\otimes h)=0$ if $u\otimes g$ and $v\otimes h$ belong to the factors in the decomposition,
 corresponding to
distinct cones in $\Sigma_X$.
\end{thm}      
        
\begin{pf}  
The first factor $A_1(X)_{p,q}$ in the decomposition of $H^{p,q}(X)$ is the toric part of the cohomology.
We will prove that the rest is isomorphic to the residue part.

Let $u\in U^\sigma(X)$ and
$g\in R^\sigma_1(f)_{*\beta-\beta_0+\beta_1^\sigma}$.
By the definition of $U^\sigma(X)$, we can assume that
$u=b\cdot\prod_{\rho_k\subset\gamma}i^*D_k$ for some $b\in U^\sigma(X)$ and 
$\gamma\in\Sigma$
such that $\tilde\pi(\inte\gamma_k)\subset\inte\sigma$.
The map assigning $b{\varphi_{\gamma}}_!\res_\gamma(g)$ to  $u\otimes g$ is well defined,
because
$$\frac{{\varphi_{\gamma}}_!\res_{\gamma}(A)}{\mult(\gamma)}
=\prod_{\rho_k\subset\gamma'\setminus\gamma}i^*D_k{\varphi_{\gamma'}}_!\res_{\gamma'}(A)$$
for any $\gamma'\subset\gamma$ such that $\tilde\pi(\inte\gamma')\subset\inte\sigma$
and
$\tilde\pi(\inte\gamma)\subset\inte\sigma$.
To show this use the Gysin projection formula as in the proof of Proposition~\ref{p:produ}.
We omit these details.
Hence, by Proposition~\ref{p:pdes}, we have a surjective map
\begin{equation}\label{e:map}
\bigoplus_{\sigma\in\Sigma_X}
U^\sigma(X)_{s,s}
\otimes R^\sigma_1(f)_{(q-s+1)\beta-\beta_0+\beta_1^\sigma}@>\psi>>H^{p,q}_{\rm res}(X),
\end{equation}
where $s=(p+q-i-\dim\sigma+1)/2$.
To find the kernel of this map we use the Poincar\'e duality.
This duality and Proposition~\ref{p:produ}{\rm(iii)} show that the kernel cannot intersect the sum of two or more
factors on the left side of (\ref{e:map}). Therefore, we can assume that a given element from the kernel belongs to
one of the factors
$U^\sigma(X)_{s,s}\otimes R^\sigma_1(f)_{(q-s+1)\beta-\beta_0+\beta_1^\sigma}$ for some $\sigma\in\Sigma_X$.
Let $\sum_{k=1}^l u_k\otimes g_k$ be this element, where $g_1,\dots,g_l$ is a basis of 
$R^\sigma_1(f)_{(q-s+1)\beta-\beta_0+\beta_1^\sigma}$, and $u_k\in U^\sigma(X)_{s,s}$.
The composition 
$$R^\sigma_1(f)_{(q-s+1)\beta-\beta_0+\beta_1^\sigma}\otimes R^\sigma_1(f)_{(i(\sigma)-q+s)\beta-\beta_0+\beta_1^\sigma}
@>>>R^\sigma_1(f)_{(i(\sigma)+1)\beta-2\beta_0+2\beta_1^\sigma}@>\res^\sigma>>{\Bbb C},
$$
where the first arrow is a multiplication, is a nondegenerate pairing because of Proposition~\ref{p:iso}(ii),
\cite[Theorem~4.4]{m1} and the definition of  $\res^\sigma$.
Hence, there is a basis $g_1',\dots,g_l'$ of $R^\sigma_1(f)_{(i(\sigma)-q+s)\beta-\beta_0+\beta_1^\sigma}$,
 orthogonal to $g_1,\dots,g_l$ with respect to this pairing. 
The fact that $\psi(\sum_{k=1}^l u_k\otimes g_k)=0$  is equivalent to $\int_X\psi(\sum_{k=1}^l u_k\otimes g_k)\cup t=0$
for all $t\in H^*(X)$.  This is true for $t\in H_{\rm toric}^*(X)$, and, by  Proposition~\ref{p:produ}{\rm(iii)},
we can consider only $t=\psi(v\otimes g_k')$ with $v\in U^\sigma(X)_{d-i(\sigma)-s,d-i(\sigma)-s}$ and $k=1,\dots,l$.
Applying the parts (i) and (ii) of Proposition~\ref{p:produ}, we get
$$\int_X\psi(\sum_{k=1}^l u_k\otimes g_k)\cup \psi(v\otimes g_k')=\int_\ps c_{q-s}^\sigma\res^\sigma(g_kg_k')X^{i(\sigma)}u_kv.$$
Hence, $\psi(\sum_{k=1}^l u_k\otimes g_k)=0$ if and only if 
$X^{i(\sigma)}u_kv\in(P(\Sigma)+SR(\Sigma))$ for all $k$ and $v\in U^\sigma(X)_{d-i(\sigma)-s,d-i(\sigma)-s}$.
In the last condition we can take all $v\in U^\sigma(X)$, because 
 the cohomology  
$H^*(\ps)\cong{\Bbb C}[D_1,\dots,D_n]/(P(\Sigma)+SR(\Sigma))$ of the complete simplicial toric variety $\ps$
satisfies the Poincar\'e duality. 
Thus, the decomposition of $H^{*,*}(X)$ is proved, and the product structure on it follows from Proposition~\ref{p:produ}.
\end{pf}

\section{The Picard group of semiample  hypersurfaces.}\label{s:pic}
        
In this section we will explicitly describe the Picard group of semiample nondegenerate
 hypersurfaces in a complete simplicial toric variety.
Our approach  provides a better description
than the one in \cite{r} with complicated restrictions.

First, we have the following property for semiample nondegenerate
 hypersurfaces.

\begin{lem} Let $X$ be an $i$-semiample  nondegenerate hypersurface in a complete simplicial toric variety $\ps$.
Then 
$$H^k(X,{\cal O}_X)=0\quad\text{ for  }k\ne0,i-1.$$
\end{lem}

\begin{pf} 
Since 
$$H^k(X,{\cal O}_X)\cong H^{0,k}(X)\cong H_{\rm toric}^{0,k}(X)\oplus H_{\rm res}^{0,k}(X)$$ 
and $H^{0,k}(\ps)$ vanishes for
$k\ne0$, the statement is implied  by Proposition~\ref{p:pdec}.
\end{pf} 

By this lemma, 
for an $i$-semiample  nondegenerate hypersurface $X$   with $i>3$, 
 we have $$H^1(X,{\cal O}_X)=H^2(X,{\cal O}_X)=0.$$ 
Combining this with the standard exponential sequence
$$0@>>>{\Bbb Z}@>>>{\cal O}_X@>{\rm exp}>>{\cal O}^*_X@>>>0$$
is enough to conclude
$$\pic(X)_{\Bbb C}\cong H^2(X,{\Bbb C})\cong H^{1,1}(X).$$

\begin{thm}\label{t:pic} Let $X\subset\ps$ be an $i$-semiample nondegenerate hypersurface 
defined by $f\in S_\beta$, and let $i>3$.
Then 
\begin{multline*}
\pic(X)_{\Bbb C}\cong H^2(X,{\Bbb C})\cong H^{1,1}(X)\\
\cong\biggl(\bigoplus_{k=1}^n{\Bbb C}D_k\biggr)/C\bigoplus\biggl(\bigoplus_{\sigma\in\Sigma_X(i-1)}\,
\bigoplus_{\tilde\pi(\inte\rho_k)\subset\inte\sigma}{\Bbb C}D_k\otimes
R^\sigma_1(f)_{\beta-\beta_0+\beta_1^\sigma}\biggr),
\end{multline*}
where 
$C:={\rm span}_{\Bbb C}\bigl\{\sum_{k=1}^n\langle m,e_k\rangle D_k,
D_l:\, m\in M,\,\tilde\pi(\inte\rho_l)\subset\inte\sigma,\,\sigma\in\Sigma_X(i)\bigr\}$.
\end{thm}      
        
\begin{pf}  By Theorem~\ref{t:main},
$$H^{1,1}(X)\cong A_1(X)_{1,1}\oplus\biggl(\bigoplus_{\sigma\in\Sigma_X(i-1)}
A_1^\sigma(X)_{1,1}
\otimes R^\sigma_1(f)_{\beta-\beta_0+\beta_1^\sigma}\biggr).$$
One may try to use only the definitions of $A_1(X)$ and $A_1^\sigma(X)$ to deduce the result. However, we
will use some additional information to alleviate this task.
By definition, $A_1(X)_{1,1}$ is a quotient  of $H^{1,1}(\ps)\cong \pic(\ps)_{\Bbb C}$, the description
of which we know, for instance, from \cite{f1}. Also,  Proposition~\ref{p:regh} or Lemma~\ref{l:int} imply that
$X\cdot D_l=0$ in $H^*(\ps)$ for $\tilde\pi(\inte\rho_l)\subset\inte\sigma$, $\sigma\in\Sigma_X(i)$.
This shows that $H^{1,1}(X)$ is at least a quotient space of the given answer.
To prove this answer it now suffices to verify that $h^{1,1}(X)$ coincides with 
$$n-d-\sum_{\sigma\in\Sigma_X(i)}a_1(\sigma)+\sum_{\sigma\in\Sigma_X(i-1)}a_1(\sigma)
\dim R^\sigma_1(f)_{\beta-\beta_0+\beta_1^\sigma},$$
where $a_1(\sigma)$ is the number of the 1-dimensional cones $\rho_l$ such that  $\tilde\pi(\inte\rho_l)\subset\inte\sigma$.
As in \cite[Section~4]{m1}, the Gysin spectral sequence gives  the exact sequence:
$$
0@>>>H^{1,1}(H^1(X\cap\ttt))@>>>\bigoplus_{k=1}^n H^{0,0}(X\cap D_k)@>>>H^{1,1}(X)@>>>H^{1,1}(H^2(X\cap\ttt))@>>>0. 
$$
Since $i>3$, 
$H^{1,1}(H^2(X\cap\ttt))$ vanishes as in the proof of Proposition~\ref{p:pdec}(i).
On the other side, Theorem~\ref{t:lef} gives $h^{1,1}(H^1(X\cap\ttt))=h^{1,1}(H^1(\ttt))=d$.
We also have 
$$H^{0}(X\cap D_k)=H_{\rm toric}^{0}(X\cap D_k)\oplus H^0_{\rm res}(X\cap D_k).$$
The dimension of $H_{\rm toric}^{0}(X\cap D_k)$ is clearly 1, when  $X\cap D_k$ is nonempty. By Proposition~\ref{p:regh},
the last condition
is the same as  $\tilde\pi(\inte\rho_k)\subset\inte\sigma$ for some $\sigma\in\Sigma_X$ of dimension less than $i$.
Finally, Proposition~\ref{p:pdec} gives  $H^0_{\rm res}(X\cap D_k)\cong R^\sigma_1(f)_{\beta-\beta_0+\beta_1^\sigma}$,
and a simple count verifies the dimension of $H^{1,1}(X)$.
\end{pf}
        
 It would actually be  convenient for us to use a more geometric description of the Picard group of $X$.
The factor ${\Bbb C}D_k\otimes R^\sigma_1(f)_{\beta-\beta_0+\beta_1^\sigma}$ in the above theorem is isomorphic
to $H^0_{\rm res}(X\cap D_k)$, while there is  an isomorphism 
$$\pi_{\rho_k}^*:H^0_{\rm res}(Y\cap V(\sigma))\cong H^0_{\rm res}(X\cap D_k).$$
The intersection $Y\cap V(\sigma)$ is a finite set of points $p_1,\dots,p_{\vol(\Gamma_\sigma)}$ inside the 1-dimensional toric variety 
$V(\sigma)\cong{\Bbb P}^1$. Hence, 
$$H^0(Y\cap V(\sigma))\cong \bigoplus_{l=1}^{\vol(\Gamma_\sigma)}{\Bbb C}p_l.$$
It follows from the definition of the residue part that 
$H^0_{\rm res}(Y\cap V(\sigma))$ consists of  $\sum_{l=1}^{\vol(\Gamma_\sigma)}a_lp_l$ such that
$\sum_{l=1}^{\vol(\Gamma_\sigma)}a_l=0$.
Correspondingly, we get that $X\cap D_k=\bigcup_{l=1}^{\vol(\Gamma_\sigma)}X_{k,l}$ is a disjoint union of 
$\vol(\Gamma_\sigma)$ irreducible components, and 
$$H^0_{\rm res}(X\cap D_k)\cong\biggl\{\sum_{l=1}^{\vol(\Gamma_\sigma)}a_lX_{k,l}:\,
\sum_{l=1}^{\vol(\Gamma_\sigma)}a_l=0\biggr\}.$$

\begin{cor}\label{c:ald} Let $X\subset\ps$ be an $i$-semiample nondegenerate hypersurface with $i>3$.
Then 
$$
\pic(X)_{\Bbb C}
\cong\biggl(\bigoplus_{k=1}^n{\Bbb C}D_k\biggr)/C\bigoplus\biggl(\bigoplus_{\begin{Sb}\sigma\in\Sigma_X(i-1)\\
\tilde\pi(\inte\rho_k)\subset\inte\sigma\end{Sb}}\biggl\{\sum_{l=1}^{\vol(\Gamma_\sigma)}a_lX_{k,l}:\,
\sum_{l=1}^{\vol(\Gamma_\sigma)}a_l=0\biggr\}\biggr),
$$
where 
$C$ is as in Theorem~\ref{t:pic}, and $X_{k,l}$ is a connected component of $X\cap D_k$.
\end{cor}

To connect this description of the Picard group to the one in  Theorem~\ref{t:pic}, we will show which elements in the latter
correspond
to $-X_{k,l}+X_{k,l+1}$, $l=1,\dots,\vol(\Gamma_\sigma)-1$, in the former. Why we choose such a basis will be clear when
we construct a mirror map in Section~\ref{s:gmd}.

Let $X$ be determined by $f\in S_\beta$ be linearly equivalent to
a torus invariant divisor $\sum_{k=1}^nb_kD_k$ with the associated polytope 
$\Delta_D=\{m\in M_{\Bbb R}:\langle m,e_k\rangle\geq-b_k\text{ for all } k\}$.
Write $f=\sum_{m\in\Delta\cap M}a_m x^{D(m)}$ where  $x^{D(m)}=\prod_{k=1}^nx_k^{b_k+\langle m,e_k\rangle}$.
The restriction of this polynomial  $f|_{D_k}=\sum_{m\in\Gamma_\sigma\cap M}a_m x^{D(m)}$ determines
$X\cap D_k$. Since $m_1^\sigma$ is a basis element for the lattice $M_X\cap\sigma^\perp$,
the lattice points $m\in\Gamma_\sigma$ are $m_0+sm_1^\sigma$, $s=0,\dots,\vol(\Gamma_\sigma)$, where
$m_0$ is  some of the vertices of $\Gamma_\sigma$.
Then the polynomial $f|_{D_k}$ can be factored as 
$$a_{m_0+\vol(\Gamma_\sigma)m_1^\sigma}x^{D(m_0)}\prod_{s=1}^{\vol(\Gamma_\sigma)}
(\prod_{k=1}^nx_k^{\langle m_1^\sigma,e_k\rangle}-\lambda_s),$$
where $\prod_{k=1}^nx_k^{\langle m_1^\sigma,e_k\rangle}$ is a coordinate on the torus $\ttt_{\rho_k}$,
which is equal to $\lambda_s$ precisely on the connected component $X_{k,s}$ of $X\cap D_k$.
This can be seen from the description of the map $\pi$ in Proposition~\ref{p:regh}.
 
\begin{lem}\label{l:cor}
Let $X\subset\ps$ be an $i$-semiample nondegenerate hypersurface defined by $f\in S_\beta$.
Then  $-X_{k,l}+X_{k,l+1}$, $l=1,\dots,\vol(\Gamma_\sigma)-1$, corresponds to
$D_k\otimes g_l$, where
$$
g_l=\frac{(\lambda_{l+1}-\lambda_{l})f\prod_{k=1}^nx_k^{\langle m_1^\sigma,e_k\rangle}}
{(\prod_{k=1}^nx_k^{\langle m_1^\sigma,e_k\rangle}-\lambda_l)
(\prod_{k=1}^nx_k^{\langle m_1^\sigma,e_k\rangle}-\lambda_{l+1})\prod_{\tilde\pi(\rho_k)\not\subset\sigma}x_k}
\in R^\sigma_1(f)_{\beta-\beta_0+\beta_1^\sigma}.
$$
\end{lem}

\begin{pf}
It suffices to show that the class of $-X_{k,l}+X_{k,l+1}$ in $H^{0,0}(X\cap D_k)$ is represented by 
the cocycle
$${\fracwithdelims\{\}{(\lambda_{l+1}-\lambda_l)f\prod_{k=1}^nx_k^{\langle m_1^\sigma,e_k\rangle}
\langle m_1^\sigma,e_{i_0}\rangle}
{(\prod_{k=1}^nx_k^{\langle m_1^\sigma,e_k\rangle}-\lambda_l)
(\prod_{k=1}^nx_k^{\langle m_1^\sigma,e_k\rangle}-\lambda_{l+1})x_{i_0}f_{i_0}}}_{i_0},
$$
which also represents $\res_{\rho_k}(g_l)$, by Proposition~\ref{p:coc}.
Note that   
\begin{multline*}
x_{i_0}f_{i_0}|_{D_k}=a_{m_0+\vol(\Gamma_\sigma)m_1^\sigma}
x^{D(m_0)}\biggl((b_{i_0}+\langle m_{0},e_{i_0})\prod_{s=1}^{\vol(\Gamma_\sigma)}
(\prod_{k=1}^n x_k^{\langle m_1^\sigma,e_k\rangle}-\lambda_s)
\\
+
\langle m_1^\sigma,e_{i_0}\rangle\prod_{k=1}^nx_k^{\langle m_1^\sigma,e_k\rangle}\sum_{j=1}^{\vol(\Gamma_\sigma)}
\prod_{s\ne j}
(\prod_{k=1}^n x_k^{\langle m_1^\sigma,e_k\rangle}-\lambda_s)\biggr).
\end{multline*}
Since $\prod_{k=1}^n x_k^{\langle m_1^\sigma,e_k\rangle}$ has the value $\lambda_s$ on $X_{k,s}$,
the above cocycle is equal to $-1$ on $X_{k,l}$, 1 on $X_{k,l+1}$ and 0 on all other components.
\end{pf}
        
In Section~\ref{s:gmd}, we will  describe the product structure on the ring generated by $\pic(X)_{\Bbb C}$.

\section{The B-model correlation functions.}\label{s:bcor}

In \cite{m2}, we calculated a subring of the B-model chiral ring $H^*(X,\wedge^*{\cal T}_X)$, 
which contains the space $H^1(X,{\cal T}_X)$, for semiample anticanonical
nondegenerate hypersurfaces. 
Here, we notice that the ring structure for minimal hypersurfaces is related to a product of 
the roots of $A$-type Lie algebra. Let us note that a semiample anticanonical nondegenerate hypersurface $X\subset\ps$  is big and Calabi-Yau.
We  will use the products on the  chiral ring in the next section 
to propose a mirror map between the space $H^1(X,{\cal T}_X)$ and 
the Picard group of the mirror symmetric Calabi-Yau hypersurface $X^\circ$ in the Batyrev construction.

We have the  birational contraction $\pi:\ps@>>>\psx$. In \cite[Section~4]{m2}, 
we  ordered the cones $\rho_i$ inside $\sigma\in\Sigma_X(2)$,
according the way they lie in  $\sigma$: $\rho_{l_0},\dots,\rho_{l_{n(\sigma)+1}}$,
so that $\rho_{l_0}$ and $\rho_{l_{n(\sigma)+1}}$ are the edges of $\sigma$.
The corresponding divisors $D_{l_k}$ in $\ps$ intersect only if they are next to each other in the order.
Using this order we can state \cite[Theorem~7.1]{m2}.
 
\begin{thm}\label{t:chiral} \cite{m2}
Let $X\subset\ps$ be  a  semiample anticanonical nondegenerate hypersurface
defined by $f\in S_\beta$.  
Then there is a natural inclusion
$$\gamma_{\_}\oplus(\oplus\gamma^{\sigma,k}_{\_}):
R_1(f)_{*\beta}
\oplus\Biggl(\bigoplus\begin{Sb}\sigma
\in\Sigma_X(2)\end{Sb}\bigl(R^\sigma_1(f)_{(*-1)\beta+\beta_1^\sigma}\bigr)
^{n(\sigma)}\Biggr)\hookrightarrow H^*(X,\wedge^*{\cal T}_X),$$
where  the sum $\oplus\gamma^{\sigma,k}_{\_}$ is over $\rho_{l_k}$, $k=1,\dots,n(\sigma)$. Also, 
$R^\sigma_1(f)_{(q-1)\beta+\beta_1^\sigma}=0$ for $q=0,d-1$.
\end{thm}

To simplify our further calculations, in this section we make an additional assumption on $\ps$:
\begin{equation}\label{e:assu}
\mult(\sigma')=1\text{ for all } \sigma'\in\Sigma(2)\text{  such that }
\sigma'\subset\sigma\in\Sigma_X(2).
\end{equation}
 In particular, this holds if $X\subset\ps$ is a minimal Calabi-Yau hypersurface in the mirror construction of \cite{b1}.
There is a nice way to state the product structure on the subspace of $H^*(X,\wedge^*{\cal T}_X)$ in the above theorem.
For this we need to use the  $(n(\sigma)+1)$-dimensional vectors $\alpha_k=(0,\dots,0,-1,1,0,\dots,0)$ with $-1$ on the $k$-th place,
which form a basis of the root system of $A_{n(\sigma)}$-type Lie algebra. 
Also, define the following product of an arbitrary number of vectors $v_l=(v_{l,1},\dots,v_{l,n(\sigma)+1})\in{\Bbb R}^{n(\sigma)+1}$,
$l=1,\dots, L$, by the rule:
$$v_1\cdots v_L=\sum_{k=1}^{n(\sigma)+1} v_{1,k}\cdots v_{L,k}.$$

As a corollary of \cite[Theorem~6.3]{m2}, we have:

\begin{thm}\label{t:chirp}
Let $X\subset\ps$ be  a  semiample anticanonical nondegenerate hypersurface
defined by $f\in S_\beta$, and assume (\ref{e:assu}) holds.  
Then, under the identifications of Theorem~\ref{t:chiral}, we have

{\rm(i)} $\gamma_A\cup\gamma_B=\gamma_{AB}$,

{\rm(ii)}  $\gamma_{A}\cup\gamma^{\sigma,k}_{B}=\gamma^{\sigma,k}_{AB}$,

{\rm(iii)} $\gamma^{\sigma_1,k_1}_{A}\cup\gamma^{\sigma_2,k_2}_{B}=0$ if $\sigma_1\ne\sigma_2$,

{\rm(iv)} for $A,B\in R^\sigma_1(f)_{(*-1)\beta+\beta_1^\sigma}$, 
such that $AB\in R^\sigma_1(f)_{(d-3)\beta+2\beta_1^\sigma}$,
$$\gamma^{\sigma,k_1}_{A}\cup\gamma^{\sigma,k_2}_{B}=\alpha_{k_1}\cdot\alpha_{k_2}
\gamma_{\mu^{-1}(ABG^\sigma(f))}\in H^{d-1}(X,\wedge^{d-1}{\cal T}_X),$$
where $\mu:R_1(f)_{(d-1)\beta}@>>>R_0(f)_{d\beta}$ is an isomorphism defined by multiplication with $\prod_{k=1}^n x_k$,
and where 
$G^\sigma(f):=\frac{x_{l_0}f_{l_0}x_{l_{n(\sigma)+1}}f_{l_{n(\sigma)+1}}\prod_{\rho_k\not\subset\sigma} x_k}
{\mult(\sigma)\prod_{\rho_k\subset\sigma} x_k}
\in S_{3\beta-2\beta_1^\sigma}$,

{\rm(v)} for  $A,B\in R^\sigma_1(f)_{(*-1)\beta+\beta_1^\sigma}$, 
such that $ABC\in R^\sigma_1(f)_{(d-4)\beta+3\beta_1^\sigma}$,
$$\gamma^{\sigma,k_1}_{A}\cup\gamma^{\sigma,k_2}_{B}\cup\gamma^{\sigma,k_3}_{C}=
\alpha_{k_1}\cdot\alpha_{k_2}\cdot\alpha_{k_3}
\gamma_{\mu^{-1}(ABCH^\sigma(f)G^\sigma(f))}\in H^{d-1}(X,\wedge^{d-1}{\cal T}_X),$$
where 
$$H^\sigma(f):=\sqrt{-1}\sum_{m\in\sigma^\perp\cap\Delta\cap M}a_m\frac{x^{D(m)}}{\prod_{\rho_k\subset\sigma} x_k}
\in S_{\beta-\beta_1^\sigma}$$
with $\sigma\in\Sigma_X(2)$, $f=\sum_{m\in\Delta\cap M}a_m x^{D(m)}$,  and 
the polytope $\Delta$ corresponding to the anticanonical divisor $\sum_{k=1}^nD_k$.
\end{thm}

We can now see a peculiar product structure on a part of 
the B-model chiral ring $\oplus_{p,q} H^p(X,\wedge^q{\cal T}_{X})$. This ring contains modules
$\oplus_q R^\sigma_1(f)_{(q-1)\beta+\beta_1^\sigma}$, corresponding to $e_k\in\inte\sigma$,
 over the subring $\oplus_q R_1(f)_{q\beta}$ called the polynomial
part in \cite{m2}. The product structure on the sum of these modules is governed by a product of the roots
of    $A_{n(\sigma)}$-type Lie algebra. While this  appears to be surprising,
 one likely
reason for such a mystery comes from the fact that the contraction $Y=\pi(X)$ has $A_{n(\sigma)}$-type
singularities along a codimension 2 subvariety $Y\cap V(\sigma)$. The generators $e_k\in\inte\sigma$
correspond to the exceptional divisors $D_k$ of the blow up $\pi$, whose intersection matrix forms the
Dynkin diagram of $A_{n(\sigma)}$, as it was observed in \cite{kmp} for 
the 3-dimensional Calabi-Yau hypersurfaces.

While it is difficult to compute all of    
the $(d-1)$-point functions of the Calabi-Yau hypersurface, the following products follow from
Theorem~\ref{t:chirp}.

\begin{cor}\label{c:pro}
Let $X\subset\ps$ be  a  semiample anticanonical nondegenerate hypersurface
defined by $f\in S_\beta$, and assume (\ref{e:assu}) holds.  
Given $A_1,\dots,A_s\in R_1(f)_\beta$ and $B_1,\dots,B_t\in R^\sigma_1(f)_{\beta_1^\sigma}$, $s+t=d-1$,
$t\le3$, then
$$
\gamma_{A_1}\cup\cdots\cup\gamma_{A_s}\cup\gamma^{\sigma,k_1}_{B_1}
\cup\cdots\cup\gamma^{\sigma,k_{t}}_{B_{t}}=
\alpha_{k_1}\cdots\alpha_{k_t}
\gamma_{\mu^{-1}(A_1\cdots A_sB_1\cdots B_tH^\sigma(f)^{t-2}G^\sigma(f))}.
$$
\end{cor}

This corollary together with the products on the polynomial part $\oplus_q R_1(f)_{q\beta}$
of the B-model chiral ring, gives all information about the 
non-normolized $3$-point functions (Yukawa couplings)
for the  Calabi-Yau 3-fold hypersurfaces.

\section{A generalization of the monomial-divisor mirror map.}\label{s:gmd}
 
The monomial-divisor mirror map was proposed in \cite{agm}. This map was conjectured to
be the derivative of the mirror map at large radius limit points (maximally unipotent boundary points) 
between the complex and K\"ahler moduli spaces of a mirror pair $(X,X^\circ)$ of Calabi-Yau hypersurfaces
in toric varieties. Here we extend this map to the whole $\pic(X)_{\Bbb C}$ and 
$H^1(X^\circ,{\cal T}_{X^\circ})$.
Our map is supported by a compatibility of some limiting products of the chiral rings.

While Theorem~\ref{t:main} already has a description of the products on the 
$\pic(X)_{\Bbb C}$,
it is very technical to compute  $\res^\sigma$ and $p_\sigma$  directly. That is why we found 
an alternative description
of the Picard group of $X$ in Corollary~\ref{c:ald}. Using this description we will  prove the following.

\begin{pr}\label{p:sam} Let $X\subset\ps$ be an $i$-semiample nondegenerate hypersurface 
defined by $f\in S_\beta$. Given $a_1,\dots,a_s\in H^{1,1}_{\rm toric}(X)$, then 
$$a_1\cup\cdots\cup a_s\cup (D_{i_1}\otimes g_{k_1})\cup\cdots\cup(D_{i_t}\otimes g_{k_t})=
\alpha_{k_1}\cdots\alpha_{k_t}a_1\cdots a_s\cdot D_{i_1}\cdots D_{i_t},
$$
where $s+t=d-1$.
\end{pr}

\begin{pf} 
By Lemma~\ref{l:cor}, $D_{i_j}\otimes g_{k_j}$ corresponds to $-X_{i_j,k_j}+X_{i_j,k_j+1}$, 
which is the same as ${\varphi_{\rho_{i_j}}}_!\pi_{\rho_{i_j}}^*(-p_{k_j}+p_{k_j+1})$. Let us remind that
$$Y\cap V(\sigma)=\pi(X)\cap V(\sigma)=\bigcup_{k=1}^{\vol(\Gamma_\sigma)}{p_k}.$$
The product structure on 
$$H^0(Y\cap V(\sigma))\cong \bigoplus_{k=1}^{\vol(\Gamma_\sigma)}{\Bbb C}p_l$$
is an obvious one: multiplication of the corresponding components of the vectors.
It is not difficult to see that this product can be written as 
$$(-p_{k_1}+p_{k_1+1})\cup(-p_{k_2}+p_{k_2+1})=\frac{\alpha_{k_1}\cdot\alpha_{k_2}}{n(\sigma)+1}
\sum_{k=1}^{\vol(\Gamma_\sigma)}{p_k}+
\sum_{k=1}^{n(\sigma)}c_k(-p_k+p_{k+1})
$$
where the coefficients $c_k$ are determined by the decomposition
$$(\alpha_{k_1,1}\alpha_{k_2,1},\dots,\alpha_{k_1,n(\sigma)+1}\alpha_{k_2,n(\sigma)+1})=
\frac{\alpha_{k_1}\cdot\alpha_{k_2}}{n(\sigma)+1}(1,\dots,1)+\sum_{k=1}^{n(\sigma)}c_k\alpha_k.$$
Hence,
\begin{multline*}
(D_{i_1}\otimes g_{k_1})\cup(D_{i_2}\otimes g_{k_2})=
{\varphi_{\rho_{i_1}}}_!\pi_{\rho_{i_1}}^*(-p_{k_1}+p_{k_1+1})
\cup{\varphi_{\rho_{i_2}}}_!\pi_{\rho_{i_2}}^*(-p_{k_2}+p_{k_2+1})
\\
=i^*D_{i_2}\cup{\varphi_{\rho_{i_1}}}_!\pi_{\rho_{i_1}}^*((-p_{k_1}+p_{k_1+1})\cup(-p_{k_2}+p_{k_2+1}))
\\
=\frac{\alpha_{k_1}\cdot\alpha_{k_2}}{n(\sigma)+1}D_{i_1}D_{i_2}+
\sum_{k=1}^{n(\sigma)}c_k(D_{i_1}D_{i_2}\otimes g_{k}),
\end{multline*}
where we also use the same technique as in the proof of Proposition~\ref{p:produ}.

Now we can prove the proposition by induction on $t$. It is true for $t=1$: the left side
of the formula vanishes because the residue part is orthogonal to the toric part,
while the other side also vanishes by our definition of the products of the vectors. 
Assuming that the statement holds for $t-1$,
we have
\begin{multline*}
a_1\cup\cdots\cup a_s\cup (D_{i_1}\otimes g_{k_1})\cup\cdots\cup(D_{i_t}\otimes g_{k_t})
\\
=
a_1\cup\cdots\cup a_s\cup\biggl(\frac{\alpha_{k_1}\cdot\alpha_{k_2}}{n(\sigma)+1}D_{i_1}D_{i_2}+
\sum_{k=1}^{n(\sigma)}c_k(D_{i_1}D_{i_2}\otimes g_{k})\biggr)\cup
(D_{i_3}\otimes g_{k_3})\cup\cdots\cup(D_{i_t}\otimes g_{k_t})
\\
=
\frac{\alpha_{k_1}\cdot\alpha_{k_2}}{n(\sigma)+1}a_1\cdots a_s D_{i_1}D_{i_2}\cdots D_{i_t}
+\sum_{k=1}^{n(\sigma)}c_k\alpha_{k}\cdot\alpha_{k_3}\cdots\alpha_{k_t}
a_1\cdots a_s D_{i_1}D_{i_2}\cdots D_{i_t}
\\
=
\alpha_{k_1}\cdots\alpha_{k_t}a_1\cdots a_s\cdot D_{i_1}\cdots D_{i_t}.
\end{multline*}
\end{pf}

We can see a similarity between the products in the above proposition and Corollary~\ref{c:pro}.
Based on this we will propose a mirror map between $\pic(X)_{\Bbb C}$ and 
$H^1(X^\circ,{\cal T}_{X^\circ})$ for a mirror  pair $(X,X^\circ)$ of Calabi-Yau hypersurfaces.

Let us  recall the mirror symmetry construction from \cite{b1}.
We will 
describe it starting with
a  semiample nondegenerate Calabi-Yau hypersurface $X$ in a complete simplicial toric variety $\ps$.
Such a hypersurface has the anticanonical degree $[\sum_{i=1}^n D_n]$.
Theorem~\ref{t:fun}  gives a unique birational toric morphism $\pi:\ps@>>>\psx$
such that the fan $\Sigma$ is a subdivision of the fan $\Sigma_X$, $\pi_*[X]$ is ample 
and $\pi^*\pi_*[X]=[X]$.
 Moreover, $\Sigma_X$ is the normal fan of the polytope $\Delta$
associated with a torus invariant  divisor equivalent to $X$.
The push-forward $\pi_*[X]$ is again anticanonical, whence the toric variety $\psx$ is
Fano associated with the polytope
$\Delta\subset M_{\Bbb R}$ of the anticanonical divisor $\sum_{i=1}^n D_n$ on $\ps$, by Lemma 3.5.2 
in \cite{ck}.
Then, Proposition \ref{p:regh} shows that the image $Y:=\pi(X)$ is an ample nondegenerate hypersurface 
in $\psx=\pp_\Delta$.
The fact that  $\pp_\Delta$ is Fano means by Proposition 3.5.5 in \cite{ck} that the polytope 
$\Delta$ is reflexive,
i.e., its dual 
$$\Delta^\circ=\{n\in N_{\Bbb R}: \langle m,n\rangle\ge-1\text{ for } m\in\Delta\}$$
has all its vertices at lattice points in $N$, and the only  lattice point in the interior of $\Delta^\circ$ is the origin $0$.
This implies that the generators of the 1-dimensional cones of the fan $\Sigma$ are among 
the lattice points of $\Delta^\circ$.
If $X$ was a minimal Calabi-Yau (see \cite{ck}), then these generators coincide with the set  
$\Delta^\circ\cap N\setminus\{0\}$.
Now, consider the polytope $\Delta^\circ$ and the associated toric variety 
$\pp_{\Delta^\circ}$. Theorem 4.1.9 in \cite{b1} says that
an anticanonical nondegenerate hypersurface
$Y^\circ\subset\pp_{\Delta^\circ}$ is a Calabi-Yau variety with canonical singularities.
To obtain a mirror of $X$, it is natural to take a certain desingularization of $Y^\circ$. 
Batyrev introduced the following notion.

\begin{defn}
A projective birational morphism  $\varphi:W'@>>> W$ is called a {\it maximal projective crepant partial desingularization}
(MPCP-desingularization) of  $W$ 
if $\varphi$ is crepant ($\varphi^*(K_W)=K_{W'}$)
and $W'$ has only $\Bbb Q$-factorial terminal singularities.
\end{defn}

Take any simplicial subdivision $\Sigma^\circ$ of the normal fan of  $\Delta^\circ$ such that the generators of the
 1-dimensional cones of $\Sigma^\circ$ is the set $\Delta\cap M\setminus\{0\}$.
 Such a subdivision corresponds to
a toric birational morphism
$\pp_{\Sigma^\circ}@>>>\pp_{\Delta^\circ}$ which is a MPCP-desingularization of 
$\pp_{\Delta^\circ}$, by Lemma 4.1.2 in \cite{ck}.  
Theorem 2.2.24 of \cite{b1} says that the toric variety $\pp_{\Delta^\circ}$ admits at least one MPCP-desingularization. 
This also  induces a MPCP-desingularization $X^\circ\subset\pp_{\Sigma^\circ}$ 
of the anticanonical hypersurface $Y^\circ\subset\pp_{\Delta^\circ}$. 
The hypersurface $X^\circ$ is a minimal Calabi-Yau variety.
Batyrev showed the following formulas for some Hodge  numbers of the minimal 
Calabi-Yau hypersurfaces (MPCP-desingularizations):
$$h^{d-2,1}(X^\circ)=h^{1,1}(X)=l(\Delta^\circ)-d-1-\sum_{{\rm codim}\theta=1}
l^*(\theta)+\sum_{{\rm codim}\theta=2}
l^*(\theta)l^*(\theta^*),$$
where $\theta$ is a face of 
$\Delta^\circ$, $\theta^*$ is the corresponding dual face of the dual reflexive polyhedron
$\Delta$, and $l(\Gamma)$ (resp., $l^*(\Gamma)$) 
denotes the number of lattice (resp., interior lattice) points
in $\Gamma$. 
For Calabi-Yau 3-folds, this is called the topological mirror symmetry test.
Note also that since $X^\circ$ is Calabi-Yau, $h^{d-2,1}(X^\circ)=\dim H^1(X^\circ,{\cal T}_{X^\circ})$. 

The monomial-divisor mirror map can be easily described if we assume that
the polynomial $f\in S_\beta$ defining the hypersurface $X^\circ\subset\pp_{\Sigma^\circ}$
has the form
 $$\sum_{m\in(\Delta^\circ\cap M)_0}a_m x^{D(m)},$$
 where $(\Delta^\circ\cap M)_0$
 is the set of all lattice points in $\Delta^\circ$ which do not lie in the interior of any facet
of $\Delta^\circ$. Such a hypersurface has been called {\it simplified} in \cite{agm}, \cite{ck}.
In this case, the map between the toric part $H^{1,1}_{\rm toric}(X)$
and the polynomial part $R_1(f)_\beta$ of $H^1(X^\circ,{\cal T}_{X^\circ})$
  is defined by sending the divisor $D_k$, corresponding to the lattice point $e_k$ in 
$(\Delta^\circ\cap M)_0$, to the monomial $x^{D(e_k)}$.
One can verify that this map is well defined and induces an isomorphism.
We should note that if the hypersurface is not simplified the mirror map is complicated to describe
explicitly. 

Our extension of the  monomial-divisor mirror map is as follows.
Assign to the divisors $-X_{k,l}+X_{k,l+1}=D_k\otimes g_l\in H^{1,1}(X)$, 
for $e_k\in \inte\sigma$, $\sigma\in\Sigma_X(d-1)$ and $l=1,\dots,\vol(\Gamma_\sigma)-1$,
the elements $\gamma^{\sigma^\circ,l}_{A}$ with 
$$A=\frac{x^{D(e_k)}}{\prod_{\rho_k\not\subset\sigma} x_k}\in 
R^{\sigma^\circ}_1(f)_{\beta_1^{\sigma^\circ}},$$
where $\sigma^\circ\in\Sigma_{X^\circ}(2)$ is the cone over the 1-dimensional face $\Gamma_\sigma$.
Using the formula for the Hodge numbers and comparing the dimensions, one deduces
that $\dim R^{\sigma^\circ}_1(f)_{\beta_1^{\sigma^\circ}}$ coincides with the number of $e_k\in \inte\sigma$.
 Therefore, our map is an isomorphism.
We conjecture that this map is the derivative of the mirror map at the large radius limit points
of the
K\"ahler and complex moduli spaces
of the pair  of Calabi-Yau hypersurfaces.
One consequence of mirror symmetry is  that the derivative of the mirror map between the moduli spaces 
gives a ring isomorphism between the quantum cohomology of $X$ and the B-model chiral ring
of $X^\circ$ for a mirror pair of smooth Calabi-Yau hypersurfaces.
In particular, the compatibility of the products should also be valid at the limit points.
The limit of the product of the quantum cohomology ring is the usual cup product on the cohomology.
On the other hand,  we have the products of the B-model chiral ring in Corollary~\ref{c:pro}. 
While it is very technical to compute explicitly the
limits of  non-vanishing products, the limits of the vanishing products are zero.
The products in Corollary~\ref{c:pro} vanish when  
there are 2 non-consecutive numbers among $k_1,\dots,k_{t}$ or if  $k_1,\dots,k_{t}$ are the same number
and $t$ is odd. The same holds for the products of the cohomology ring in Proposition~\ref{p:sam}.
This verifies  the compatibility of our generalization 
of the monomial-divisor mirror map with some of the limiting products.

Since  the vanishing  products hold on the B-model chiral ring far from the limiting points in
the complex moduli space, we believe that a similar statement should be true for the quantum 
products of the mirror Calabi-Yau hypersurface. 
Thus, we conjecture that 
\begin{multline*}
a_1*\cdots* a_s* (D_{i_1}\otimes g_{k_1})*\cdots*(D_{i_t}\otimes g_{k_t})
\\=
a_1*\cdots* a_s*(-X_{i_1,k_1}+X_{i_1,k_1+1})*\cdots*(-X_{i_t,k_t}+X_{i_t,k_t+1})
=0,
\end{multline*}
where $*$ denotes the small quantum product on $H^*(X)$ (see \cite{ck}),
if 
 $|k_a-k_{b}|>1$ for some $a,b$ or if  $k_1=\cdots=k_{t}$ and $t$ is odd.

\end{document}